\documentclass[12pt,leqno,twoside]{article} 
\usepackage{bezier,eufrak}
\input{rogg_kue.sty}
\newcommand{\enger}{\setlength{\arraycolsep}{1.5pt}
                           \renewcommand{\arraystretch}{0.7} }
\newcommand{\weiter}{\setlength{\arraycolsep}{3pt}
                           \renewcommand{\arraystretch}{1.2} }
\weiter

\begin{document}
\title{Elementary divisors of Gram matrices of certain Specht modules}
\author{M.\ K\"unzer, G.\ Nebe}
\maketitle

\begin{abstract}
The elementary divisors of the Gram matrices of Specht modules $S^\lambda$ over the symmetric group are determined for two-row partitions and for two-column partitions $\lambda$. 
More precisely, the subquotients of the Jantzen filtration are calculated using Scha\-per's formula. Moreover, considering a general 
partition $\lambda$ of $n$ at a prime $p > n - \lambda_1$, the only possible non trivial composition factor of $S_{\sF_p}^\lambda$ is induced by the morphism of Carter and Payne, 
as shown by means of Kleshchev's modular branching rule. This enables the Jantzen filtration to be calculated in this case as well.
\end{abstract}

\renewcommand{\thefootnote}{\fnsymbol{footnote}}
\renewcommand{\thefootnote}{\arabic{footnote}}

\begin{footnotesize}
\renewcommand{\baselinestretch}{0.7}
\parskip0.0ex
\tableofcontents
\parskip1.2ex
\renewcommand{\baselinestretch}{1.0}
\end{footnotesize}

\setcounter{section}{-1}

\section{Introduction}

\subsection{Problem}

Specht modules $S^\lambda$ are combinatorially defined $\Z\Sl_n$-modules, indexed by partitions $\lambda$ of $n$, which yield a complete set of pairwise nonequivalent
ordinary irreducible representations of the symmetric group $\Sl_n$ after scalar extension to $\C$. A Specht module $S^\lambda$ carries a nondegenerate $\Sl_n$-invariant
bilinear form, inducing an embedding into its $\Z$-linear dual, $S^\lambda\hra S^{\lambda,\ast}$. The problem is to determine the structure of the quotient
\[
S^{\lambda,\ast}/S^\lambda
\]
as an abelian group. Reformulated, we ask for the elementary divisors of the Gram matrix of this bilinear form on $S^\lambda$, or, for short, for the elementary divisors of 
$S^\lambda$.

\subsection{Known Results}

\subsubsection{Simple modules}

Let $p$ be a prime. Denote by $\lambda'$ the transposed partition of $\lambda$. The number of elementary divisors of $S^\lambda$ not divisible by $p$ is either zero, 
or the dimension of a simple $\F_p\Sl_n$-module.

{\bf Theorem} {\bf\cite[\rm 11.5]{J}}. {\it Let $D_{\sF_p}^\lambda$ be the image of $S^\lambda\lra S^{\lambda,\ast}/pS^{\lambda,\ast}$. If $\lambda$ is
$p$-regular, that is, $\lambda'_i - \lambda'_{i+1} < p$ for all $i\geq 1$, then $D_{\sF_p}^\lambda$ is a simple $\F_p\Sl_n$-module. Up to isomorphism, all simple $\F_p\Sl_n$-modules
occur this way. If $\lambda$ is $p$-singular, we have $D_{\sF_p}^\lambda = 0$.}

Allowing for $D_{\sF_p}^\lambda$ to denote the zero module if $\lambda$ is not a $p$-regular partition turns out to be a convenient convention.

\subsubsection{Schaper's formula}

First of all, the product of the elementary divisors of $S^\lambda$ is known.

{\bf Theorem} {\bf\cite[\rm p.\ 224]{JM}}. {\it There is an explicit combinatorial formula for the determinant of the Gram matrix of $S^\lambda$.}

More precisely, given a prime $p$, the quotient $S^{\lambda,\ast}/S^\lambda$ is expressible as a linear combination of Specht modules in the Grothendieck group of $\Z\Sl_n$-modules of 
finite length by {\sc Schaper}'s formula. Given a commutative ring $A$, we abbreviate $S_A^\lambda := A\ts_\sZ S^\lambda$.

{\bf Theorem} {\bf\cite[\rm p.\ 60]{Sch}}, cf.\ {\bf\cite[\rm Cor.\ 5.33]{Mathas}}. {\it There are combinatorially determined integral coefficients $\alpha_\mu$ such that
\[
[S_{\sZ_{(p)}}^{\lambda,\ast}/S_{\sZ_{(p)}}^\lambda]\; =\; \sum_{i\geq 1} [S_{\sF_p}^\lambda(i)]\; =\; \sum_\mu \alpha_\mu [S_{\sF_p}^\mu]\;\in\; K_0(\modfr\Z_{(p)}\Sl_n)\; ,
\]
where $S_{\sF_p}^\lambda(i)$ denotes the $i$th piece of the Jantzen filtration of $S_{\sF_p}^\lambda$.}

Together with calculations of decomposition numbers due to {\sc James,} {\sc Williams,} {\sc To Law,} {\sc Benson,} {\sc M\"uller} et al. {\bf\cite{B,J,JW,Atlas,M,M18,TL}}, 
Schaper's formula represents our principal tool.

\subsubsection{Numerical results}

There is an estimate for the first elementary divisor, found by {\sc James} in the course of the construction of the simple modules $D_{\sF_p}^\mu$.

{\bf Lemma} {\bf\cite[\rm 10.4]{J}}. {\it The product $\prod_{i\geq 1} (\lambda'_i - \lambda'_{i+1})!$ divides the first elementary divisor of $S^\lambda$. In turn, the first 
elementary divisor of $S^\lambda$ divides the product $\prod_{i\geq 1} (\lambda'_i - \lambda'_{i+1})!^i$.}

For instance, $3!$ divides the first elementary divisor $12$ of $S^{(2^3)}$, which in turn divides $3!^2$. The {\it James factor} $\prod_{i\geq 1} (\lambda'_i - \lambda'_{i+1})!$ will 
reappear constantly.

Numerically, the relation between the elementary divisors of $S^\lambda$ and $S^{\lambda'}$ has been known.

{\bf Proposition} {\bf\cite[\rm 6.2.10]{K99}}. {\it Let $n_\lambda := \rk_\sZ S^\lambda$, and let $i\in [1,n_\lambda]$. The product of the $i$th elementary divisor of $S^\lambda$ and 
the $(n_\lambda+1-i)$th elementary divisor of $S^{\lambda'}$ yields $n!/n_\lambda$.} 

In particular, the elementary divisors of $S^\lambda$ and $S^{\lambda'}$ mutually determine each other. We shall give a module version of this relation in terms
of Jantzen subquotients (\ref{PropTr1}).

\subsubsection{Related work}

{\sc Grabmeier} used {\sc Schaper'}s analoguous formula for the Weyl modules over the Schur algebra as an ingredient to determine the graduated hull of $p$-adic Schur algebras
\mb{\bf\cite[\rm 11.13]{Grabm}}. 

{\sc Kleshchev} and {\sc Sheth} {\bf\cite[\rm 3.4]{KlSh}}, and independently, {\sc Reuter} {\bf\cite[\rm 4.2.22]{Reu}}, described the submodule structure of $S_{\sF_p}^{(n-m,m)}$.

\subsection{Results}

\subsubsection{Two-row partitions}

Let $n\geq 1$, let $0\leq m\leq n/2$ and let $p$ be a prime. Since the decomposition numbers of $S_{\sF_p}^{(n-m,m)}$ are in $\{ 0,1\}$ by James' formula, the Jantzen filtration
may be calculated by means of Schaper's formula.

{\bf Theorem} (\ref{CorT4}). {\it The multiplicities of the simple modules in the subquotients of the Jantzen filtration of $S_{\sF_p}^{(n-m,m)}$ are determined. In particular,
the elementary divisors of $S^{(n-m,m)}$ are calculated.}

Moreover, combining arguments of {\sc Plesken} {\bf\cite{P93}} and {\sc Wirsing} {\bf\cite{W}}, we show that if $m\geq 3$, then $S_\sQ^{(n-m,m)}$ does not contain a unimodular 
$\Z\Sl_n$-lattice, that is, a lattice $X$ satisfying $X\iso X^\ast$ (\ref{ThU3}). For $m\in \{ 1,2\}$, unimodular lattices do occur and have been classified by {\sc Plesken} 
{\bf\cite[\rm p.\ 98 and II.5]{P93}}.

\subsubsection{At a large prime}

Suppose given a partition $\lambda$ of $n$ and a prime $p > n - \lambda_1$. Using the theorem of {\sc Carter} and {\sc Payne} {\bf\cite[\rm p.\ 425]{CP}}, the direction of the
Carter conjecture proven by {\sc James} and {\sc Murphy} {\bf\cite[\rm p.\ 222]{JM}}, as well as {\sc Kleshchev}'s modular branching rule {\bf\cite[\rm 0.6]{Kl}}, the Jantzen 
filtration of $S_{\sF_p}^\lambda$ may be calculated.

{\bf Theorem} (\ref{ThL4}). {\it If $p$ does not divide a first row hook length in the range $[1,\lambda_2]$, then $S_{\sF_p}^\lambda$ is simple. If $p$ divides the first row hook 
length $h_t$ of the node $(1,t)$, $t\in [1,\lambda_2]$, then $[S_{\sF_p}^\lambda] = [D_{\sF_p}^\lambda] + [D_{\sF_p}^{\lambda[t]}]$, where $\lambda[t]$ is the partition arising from 
$\lambda$ by the according Carter-Payne box shift. The constituent $[D_{\sF_p}^{\lambda[t]}]$ lies in the $v_p(h_t)$th Jantzen subquotient.}

\subsubsection{Explicit diagonalization}

The results mentioned so far are based on Schaper's formula, so that no diagonalizing bases can be deduced. In general, an explicit diagonalization seems to be complicated.
For hook partitions, however, it is easier to diagonalize directly (\ref{ThHook}) than to apply Schaper's formula, as has already been remarked by {\sc James} and {\sc Mathas} 
[unpublished].

Moreover, for $S^{(2^2,1^{n-4})}$, we give bases essentially diagonalizing the Gram matrix (\ref{Th2^2}). In general, it might be worthwhile to employ modular morphisms 
in order to fully solve the diagonalization problem (cf.\ \ref{ConM5}). But note that from such a complete solution, bases for the simple $\F_p\Sl_n$-modules would ensue.

\subsubsection{The scope of Schaper's formula}

Some three- and four-part partitions are treated in section \ref{Sec34part}, partly conjecturally. The problem that remains to be solved, once all decomposition numbers
and Jantzen subquotients are known, is the following. 

Let $p$ be a prime, let $e^\mu$ be a primitive idempotent of $\Z_{(p)}\Sl_n$ belonging to $D_{\sF_p}^\mu$ and let 
$\eps^\lambda$ be the central-primitive idempotent of $\Q\Sl_n$ belonging to $S_\sQ^\lambda$. The investigation of the elementary divisors of $S^\lambda\hra S^{\lambda,\ast}$
can be reduced to the consideration of 
\[
S_{\sZ_{(p)}}^\lambda e^\mu\;\hra\; S_{\sZ_{(p)}}^{\lambda,\ast}e^\mu\; . 
\]
This is a $\eps^\lambda e^\mu\Z_{(p)}\Sl_n e^\mu$-linear map, the determinant of which can be deduced from Schaper's formula. So we are led to consider the $\Z_{(p)}$-order 
$\eps^\lambda e^\mu\Z_{(p)}\Sl_n e^\mu$, which is of rank \mb{$[S_{\sF_p}^\lambda : D_{\sF_p}^\mu]^2$} over $\Z_{(p)}$.

In particular, if $[S_{\sF_p}^\lambda : D_{\sF_p}^\mu] = 1$, we obtain $\eps^\lambda e^\mu\Z_{(p)}\Sl_n e^\mu\iso\Z_{(p)}$, which enables us to calculate elementary divisors.
If $[S_{\sF_p}^\lambda : D_{\sF_p}^\mu] > 1$, however, Schaper's formula alone is too coarse.

In another disguise, concerning the distribution of the multiplicity of a simple module $D_{\sF_p}^\mu$ over the Jantzen subquotients, Schaper's formula gives the value of a certain 
sum of which one needs to know the summands (cf.\ \ref{LemN1} iii). If $[S_{\sF_p}^\lambda : D_{\sF_p}^\mu] = 1$, then this sum consists of a single nonzero summand, that is,
the Jantzen subquotient $D_{\sF_p}^\mu$ appears in is determined. But this purely numerical point of view hides the role of the ring 
$\eps^\lambda e^\mu\Z_{(p)}\Sl_n e^\mu$ in case $[S_{\sF_p}^\lambda : D_{\sF_p}^\mu] > 1$.

\subsection{Acknowledgements}
We thank {\sc F.\ L\"ubeck} for access to his immense computations {\bf\cite{L}}, which served as our indispensable guideline. We thank {\sc A.\ Mathas} for help with several 
details.  We thank {\sc A.\ Kleshchev} for help with modular branching (\ref{ThL4}). We thank {\sc E.\ Wirsing} for help with Pell equations (\ref{ThU3}). We thank {\sc J.\ M\"uller}
for making {\bf\cite{M18}} available to us. The first author would like to thank the IGD Lyon for kind hospitality during the first part of this project.

\subsection{Conventions}
 
\begin{footnotesize}
\begin{itemize}
\item[(i)] Composition of maps is written on the right, $\lraa{a}\lraa{b} = \lraa{ab}$. Exception is made for `standard' maps, such as traces, characters \dots\ Unless mentioned
otherwise, a module is a finitely generated right module.
\item[(ii)] For $a,b\in\Z$, we denote by $[a,b] := \{c\in\Z\;|\; a\leq c\leq b\}$ the integral interval.
\item[(iii)] If $A$ is an assertion, which might be true or false, we let $\{ A\}$ equal $1$ if $A$ is true, and $0$ if $A$ is false. If, in a sum, a summand has a factor $\{ A\}$ 
attached, this summand is zero if $A$ is false, regardless whether it is well defined or not. 
\item[(iv)] Let $m$ be a positive integer, let $\bar a$ denote the residue class of $a$ modulo $m$. The assertion $\bar a\in \{ \bar b_1,\dots, \bar b_l\}$ will be written as 
$a\con_m b_1,\dots, b_l$. In particular, $a\con_m b$ means $a - b \in m\Z$. 
\item[(v)] Let $p$ be a prime, let $n\geq 1$. The $\F_p\Sl_n$-module $D_{\sF_p}^\mu$ is defined for any composition $\mu$ of $n$, and we let it be zero if $\mu$ is not a $p$-regular
partition. If $\mu$ is a $p$-regular partition, it is defined in {\bf\cite[\rm 11.2]{J}}. 
\item[(vi)] Given a partition $\lambda$, its transpose is denoted by $\lambda'$, i.e.\ $j\leq\lambda_i$ $\equ$ $i\leq\lambda'_j$. Given a $\lambda$-tableau $[a]$, its transpose
is a $\lambda'$-tableau denoted by $[a']$.
\item[(vii)] The sign of $\sigma\in\Sl_n$ is denoted by $\eps_\sigma$.
\item[(viii)] Let $\Jac(B)$ denote the Jacobson radical of a ring $B$.
\item[(ix)] The binomial coefficient $\smatze{a}{b}$ is zero unless $a\geq b\geq 0$.
\item[(x)] Given a ring $A$, we denote by $\modfr A$ the abelian category of $A$-modules of finite length, and by $K_0(\modfr A)$ its Grothendieck group. So $K_0(\modfr A)$
is the quotient of the free abelian group on the objects of $\modfr A$ modulo a relation $X - (X' + X'')$ for each short exact sequence $0\lra X'\lra X \lra X''\lra 0$ in $\modfr A$. 
The image of $X\in\Ob\modfr A$ in $K_0(\modfr A)$ shall be denoted by $[X]$. Cf.\ \cite[\rm 3.1.6]{Rosen}.
\end{itemize}
\end{footnotesize}

\section{General considerations}
\label{SecGen}

\subsection{Situation}

\subsubsection{An order $\Lambda$}

Let $R$ be a discrete valuation ring with maximal ideal generated by $\pi$, let $K$ be its field of fractions. Given an $R$-module $M$, we denote its reduction 
by $\b M := M/\pi M$. 

Let $A$ be a simple $K$-algebra (to which we may pass from a semisimple $K$-algebra, cf.\ \ref{RemS1}). Let $\Lambda\tm A$ be a full $R$-order in $A$, that is, 
an $R$-algebra that spans $A$ as a vectorspace such that $\Lambda$ is finitely generated free as a module over $R$. Let 
$
1_\Lambda = e_1 + \cdots + e_k
$
be an orthogonal decomposition into primitive idempotents of $\Lambda$, which correspond to indecomposable projective $\Lambda$-modules $P_j := e_j\Lambda$. 

\begin{Remark}
\label{RemS1}\rm
Suppose given a semisimple $K$-algebra $B$, a full $R$-order $\Xi$ in $B$ and a simple $\Xi$-lattice $X$, that is, a $\Xi$-module that is finitely generated 
free over $R$ such that $KX := K\ts_R X$ is a simple $B$-module. Let $\eps\in B$ be the central-primitive idempotent that acts as identity on $KX$. Then $X$ 
remains a simple lattice over the quasiblock $\Lambda = \Xi\eps$, which is a full $R$-order in $A = B\eps$, and to which we may reduce the situation for our 
purposes. Note that in general, the quasiblock $\Lambda = \Xi\eps$ will not be a direct summand of $\Xi$.

If $Q$ is an indecomposable projective $\Xi$-module, then $Q\eps\; (\tm KQ)$ is either zero or an indecomposable projective $\Lambda$-module. All indecomposable 
projective $\Lambda$-modules are obtained this way.

For instance, if $G$ is a finite group and $K$ a field of characteristic $0$, we may take $B = KG$, $\Xi = RG$ and $\eps$ the rational central-primitive 
idempotent belonging to the simple $KG$-module $KX$. 
\end{Remark}

\subsubsection{A suborder $\Delta\tm\Lambda$}

Let $P_{j_l}$ be representatives of the isomorphism classes of the indecomposable projective $\Lambda$-modules, 
$l\in [1,m]$, and let 
\[
f_l \; :=\; \sumd{j'\in [1,k],\; P_{j'}\,\iso\, P_{j_l}}e_{j'}\; . 
\]
Then the $f_l + \Jac(\Lambda)$ are the central-primitive idempotents of $\Lambda/\Jac(\Lambda)$.
Letting $D_l := e_{j_l}\Lambda/e_{j_l}\Jac(\Lambda)$, a complete system of inequivalent simple $\Lambda$-modules is given by $\{ D_l\; |\; l\in [1,m]\}$, and 
$f_l$ operates as identity on $D_l$.

This allows to introduce an $R$-suborder 
\[
\Delta \; :=\; \prod_{l\in [1,m]} f_l \Lambda f_l \;\tm\; \Lambda
\]
that has the same Grothendieck group as $\Lambda$, but that allows to decompose modules into smaller pieces. Since 
$\Lambda = \Ds_{l,l'\in [1,m]} f_l \Lambda f_{l'}$ has Jacobson radical 
\[
\Jac(\Lambda) = \left(\Ds_{l} \Jac(f_l \Lambda f_l)\right)\;\ds\;\left(\Ds_{l\neq l'} f_l \Lambda f_{l'}\right)\; ,
\]
the inclusion $\Delta\tm\Lambda$ induces an isomorphism
\[
\Delta/\Jac(\Delta)\;\lraiso\;\Lambda/\Jac(\Lambda)\; ,
\]
whence restriction yields
\[
K_0(\modfr\Delta)\;\llaiso\;K_0(\modfr\Lambda)\; .
\]

\subsubsection{Multiplicities}

Let $X$ be a simple $\Lambda$-lattice. There are various multiplicities to be considered.

By $[\b X : D_l]$ we denote the multiplicity of $D_l$ in a composition series of $\b X$ in the sense of Jordan-H\"older.
By $[\Lambda : P_j]$ we denote the multiplicity of $P_j$ as a direct summand of $\Lambda$, which is well defined since Krull-Schmidt holds for projective 
$\Lambda$-modules by Nakayama's Lemma. Finally, by $[KP_j : KX]$ we denote the multiplicity of $KX$ in $KP_j$ as a Jordan-H\"older constituent, or, 
equivalently, as a direct summand.

Denote $E := \dim_K\End_A KX$ and $E_l := \dim_{\b R} \End_\Lambda D_l$. Brauer reciprocity holds,
\[
E\cdot [KP_{j_l} : KX] \; = \; E_l\cdot [\b X : D_l]\; ,
\]
since the $K$-dimension of $KX e_{j_l}$ and the $\b R$-dimension of $\b X e_{j_l}$ coincide, and since the latter is calculated by the right hand side, as 
we see after passing to the case $\Lambda = \Delta f_l$, and subsequently reducing moritaequivalently to the case $e_{j_l} = 1$. Moreover, 
\[
E_l\cdot [\Lambda : P_{j_l}]\; =\; \dim_{\b R} D_l\; .
\]
In particular, $f_l$ is the sum of $E_l^{-1}\cdot\dim_{\b R} D_l$ primitive idempotents.

More generally, if $M$ is a $\Delta$-module of finite length over $R$, then $[M : D_l]$ denotes the Jordan-H\"older multiplicity of $D_l$ in $M$ as a 
$\Delta$-module. If $M$ happens to be the restriction of a $\Lambda$-module to $\Delta$, this is the same as the multiplicity of $D_l$ in $M$ as a 
$\Lambda$-module. So we do not need to specify whether the multiplicity refers to $\Delta$ or to $\Lambda$.
\begin{Remark}
\label{RemS2}\rm
If there exists a $l\in [1,m]$ such that $E_l = 1$ and $[\b X : D_l] = 1$, then $E = 1$. In other words, if some $D_l$ is absolutely simple and appears with 
multiplicity $1$ in $\b X$, then $KX$ is absolutely simple, too. 

Conversely, if $KX$ is absolutely simple and $[KP_{j_l} : KX] = 1$, then $D_l$ is absolutely simple.
\end{Remark}

\subsection{Jantzen filtration}

\begin{footnotesize}
\begin{quote}
We recall the basic facts concerning {\sc Jantzen}'s filtration arising from an embedding of simple $\Lambda$-lattices. This filtration is a tool to compare 
decomposition numbers with Jordan-H\"older multiplicities in the quotient of this embedding.
\end{quote}
\end{footnotesize}

Let $X\lraa{\phi} Y$ be a nonzero $\Lambda$-linear map between simple lattices $X$ and $Y$, thus necessarily injective. There exists an $N\geq 0$, which we 
choose and fix, such that $\pi^N Y\tm X\phi$.

\begin{Definition}
\label{DefH0}\rm
For $i\geq 0$, we let 
\[
\b X(i) \; :=\; ((\pi^i Y)\phi^{-1} + \pi X)/\pi X\;\tm\; \b X
\]
be the $i$th piece of the {\it Jantzen filtration} of $\b X$ with respect to $X\lraa{\phi} Y$. In particular, $\b X(0) = \b X$, and $\b X(N+1) = 0$. Given 
$l\in [1,m]$, the $i$th {\it Jantzen multiplicity} of $D_l$ is given by
\[
\theta_{l,i} \; :=\; [\b X(i)/\b X(i+1) : D_l]\; .
\]
\end{Definition}

\begin{Remark}
\label{RemN0}\rm
Note that $X(i) \iso (\pi^i Y\cap X\phi)/(\pi^i Y\cap \pi X\phi)$, and thus
\[
\frac{\b X(i)}{\b X(i+1)}\;\iso\; \frac{\pi^i Y\cap X\phi}{(\pi^{i+1}Y\cap X\phi) + (\pi^i Y\cap\pi X\phi)}\; .
\]
If we consider this subquotient as an $R$-module, by an elementary divisor decomposition we may assume that $X = Y = R$, and that $\phi$ is given by 
multiplication by $\pi^j$ for some $j\geq 0$. Then $\b X(i)/\b X(i+1)$ is isomorphic to $\b R$ if $i = j$, and to zero otherwise. Returning to the general case, 
we obtain therefore
\[
Y/X\phi\;\iso_R\; \Ds_{i\geq 0} (R/\pi^i R)^{\dim_{\b R} \b X(i)/\b X(i+1)} 
\;\iso_R\; \Ds_{l\,\in\, [1,m]} \left(\Ds_{i\geq 0} (R/\pi^i R)^{\theta_{l,i}}\right)^{\dim_{\b R} D_l}\; .
\]
\end{Remark}

\pagebreak[2]

\begin{Lemma}[Jantzen's Lemma]
\label{LemN1}
\Absit
\begin{itemize}
\item[(i)] In $K_0(\modfr\Lambda)$, we have $[Y/X\phi] = \sum_{i\geq 1} [\b X(i)] = \sum_{i\geq 0} i\, [\b X(i)/\b X(i+1)]$. 

\item[(ii)] In $K_0(\modfr\Delta)$, we have $[Y f_l/X\phi f_l] = \sum_{i\geq 1} [\b X(i) f_l] = \sum_{i\geq 0} i\, [\b X(i)f_l/\b X(i+1)f_l]$ for any 
$l\in [1,m]$. 

\item[(iii)] Given $l\in [1,m]$, we have
\[
\begin{array}{lclcl}
\sum_{i\geq 0} \theta_{l,i}  & = & [\b X : D_l]     & = & [\b X f_l : D_l] \\
\sum_{i\geq 0} i\theta_{l,i} & = & [Y/X\phi : D_l]  & = & [Y f_l/X\phi f_l : D_l]\; . \\
\end{array}
\]
\end{itemize}

\rm
Assertion (ii) follows from (i) by decomposition in $\modfr\Delta$. The second formula in (iii) follows from (i, ii), the first follows from the definition of 
$\theta_{l,i}$.

It remains to prove (i). There is a filtration of $\Lambda$-modules
\[
\pi^N X\phi \; =\; \pi^{N - 0} X\phi \cap\pi^N Y \;\tm\; \pi^{N - 1} X\phi \cap\pi^N Y  \;\tm\; \cdots \;\tm\; \pi^{N - N} X\phi \cap\pi^N Y = \pi^N Y\; ,
\]
the $i$th subquotient of which is
\[
\frac{\pi^{N - i} X\phi \cap\pi^N Y}{\pi^{N - (i-1)} X\phi \cap\pi^N Y}\;\iso\; \frac{X\phi \cap\pi^i Y}{\pi X\phi \cap\pi^i Y} \; =\; \b X(i)   
\]
for $i\in [1,N]$.
\end{Lemma}

\begin{Corollary}
\label{CorN2}
Suppose given $l\in [1,m]$. If $[\b X : D_l] = 1$, then 
\[
\theta_{l,i} \; =\;
\left\{
\begin{array}{ll}
1 & \mb{if $\;\; i = [Y/X\phi : D_l]$}\; , \\
0 & \mb{otherwise}\; . \\
\end{array}
\right. 
\]
\end{Corollary}

\begin{Corollary}
\label{CorN3}
Suppose given $l\in [1,m]$. Let $s\geq 1$. If $\theta_{l,i} = 0$ for $i\in [0,s-1]$ and if $[Y/X\phi : D_l] \leq s[\b X : D_l]$, then 
\[
\theta_{l,i} \; =\;
\left\{
\begin{array}{ll}
[\b X : D_l] & \mb{if $\;\; i = s$}\; , \\
0 & \mb{otherwise}\; . \\
\end{array}
\right. 
\]

\rm
In fact, the assumptions yield $\sum_{i\geq s+1} (i-s)\theta_{l,i} \leq 0$.
\end{Corollary}

\subsection{The reverse embedding}

Suppose given $\Lambda$-linear maps $(X\lraa{\phi} Y\lraa{\psi} X) = (X\lraa{a} X)$, where $a\in R$ and $v_\pi(a) = N$. Note that given $\phi$ and $N$ as above, 
such a map $\psi$ exists. Consider the Jantzen filtration $\b X(i)$ with respect to $\phi$ and the Jantzen filtration $\b Y(j)$ with respect to $\psi$.

\begin{Remark}
\label{RemR1}
Given $i\in [0,N]$, we have
$\;
\fracd{\b X(i)}{\b X(i+1)}\; \iso\; \fracd{\b Y(N-i)}{\b Y(N-i+1)}\; .
$
For $i > N$, we have $\b X(i)/\b X(i+1) = 0$.

\rm
In fact,
\[
\begin{array}{l}
\fracd{\b X(i)}{\b X(i+1)}\;\auf{\mb{\scr (\ref{RemN0})}}{\iso}\;\fracd{\pi^i Y\cap X\phi}{(\pi^{i+1}Y\cap X\phi) + (\pi^i Y\cap\pi X\phi)} \vspace*{1mm}\\
\hspace*{10mm}\iso\; \fracd{Y\psi\cap \pi^{N-i}X}{(\pi Y\psi\cap \pi^{N-i} X) + (Y\psi\cap\pi^{N-i+1} X)}\;\auf{\mb{\scr (\ref{RemN0})}}{\iso}\; 
\fracd{\b Y(N-i)}{\b Y(N-i+1)}\; . \\
\end{array}
\]
\end{Remark}

\subsection{Block diagonalization}

\begin{quote}
\begin{footnotesize}
We shall describe, to a certain extent, the block matrices of the block diagonalization of $X\lraa{\phi} Y$ resulting from a decomposition of $1_\Lambda$ into 
orthogonal primitive idempotents. The problem that remains, once all decomposition numbers and Jantzen multiplicities determined, is to find the elementary 
divisors of the blocks of size $> 1$. This amounts to a study of the local and rationally simple $R$-algebras $e_j\Lambda e_j$, which in concrete examples 
seems to be difficult to get a grip on. 
\end{footnotesize}
\end{quote}

\begin{Lemma}
\label{LemN5}
For $l\in [1,m]$, we have 
\[
Y e_{j_l}/X\phi e_{j_l}\;\iso_R\; \Ds_{i\geq 0} (R/\pi^i R)^{\theta_{l,i}\cdot E_l}\; .
\]

\rm
In fact,
\[
\begin{array}{rcl}
Y e_{j_l}/X\phi e_{j_l}
& \iso_R & \Ds_{i\geq 0} (R/\pi^i R)^{\dim_{\b R}\b X(i)e_{j_l}/\b X(i+1)e_{j_l}}\\
& =      & \Ds_{i\geq 0} (R/\pi^i R)^{(\dim_{\b R}\b X(i)f_l/\b X(i+1)f_l)\cdot E_l\cdot (\dim_{\b R} D_l)^{-1}}\\
& =      & \Ds_{i\geq 0} (R/\pi^i R)^{[\b X(i)f_l/\b X(i+1)f_l : D_l]\cdot E_l}\\
& =      & \Ds_{i\geq 0} (R/\pi^i R)^{\theta_{l,i}\cdot E_l}\; .\\
\end{array}
\]
\end{Lemma}

\begin{Corollary}
\label{CorN6}
Choosing $R$-linear bases in the direct summands $X e_j$ of $X$ and $Y e_j$ of $Y$, the matrix of $\phi$ appears in main diagonal block form. Given 
$l\in [1,m]$, the block belonging to $Xe_{j_l}\lraa{\phi} Ye_{j_l}$ has edge length 
\[
\dim_K KX e_{j_l} \; =\; E\cdot [KP_{j_l} : KX] \; =\; E_l \cdot [\b X : D_l]\; , 
\]
it appears with multiplicity
\[ 
[\Lambda : P_{j_l}] \; =\; E_l^{-1}\cdot\dim_{\b R} D_l
\]
and the valuation at $\pi$ of its determinant is given by
\[
E_l\cdot [Y/X\phi : D_l]\; .
\]
\end{Corollary}

\section{Elementary divisors for Specht modules}

\begin{quote}
\begin{footnotesize}
The calculation of the elementary divisors of the Gram matrix of a Specht lattice is equivalent to the determination of the quotient 
$S^{\lambda,\ast}/S^\lambda$ as an abelian group. Now, Schaper's formula expresses $S_{\sZ_{(p)}}^{\lambda,\ast}/S_{\sZ_{(p)}}^\lambda$ as a linear 
combination of Specht modules in the Grothendieck group of $\Z_{(p)}\Sl_n$-modules. Provided the decomposition numbers of the occurring Specht modules are known, 
this allows, by means of Jantzen's Lemma, to compare the decomposition of $S_{\sZ_{(p)}}^{\lambda,\ast}/S_{\sZ_{(p)}}^\lambda$ with the decomposition of 
$S_{\sF_p}^\lambda$, and thus in simple cases to determine the distribution of the simple constituents of $S_{\sF_p}^\lambda$ over the subquotients of the 
Jantzen filtration. Together with the dimensions of these simple constituents, this yields the structure of $S_{\sZ_{(p)}}^{\lambda,\ast}/S_{\sZ_{(p)}}^\lambda$ 
as an abelian group by (\ref{RemN0}).
\end{footnotesize}
\end{quote}

\subsection{Specht modules}

A {\it $\lambda$-tabloid} $\{ a\}$ is a $\lambda$-tableau $[a]$ with unordered rows. 
Let $M^\lambda$ be the free $\Z$-module on the set of $\lambda$-tabloids, carrying a structure as a $\Z\Sl_n$-module by entrywise operation of $\Sl_n$. That is, 
$M^\lambda$ is isomorphic to the permutation module on $(\Sl_{\lambda_1}\ti\cdots\ti\Sl_{\lambda_n})\< \Sl_n$. Let the $\Sl_n$-invariant bilinear form $(-,=)$ 
on $M^\lambda$ be defined by 
\[
(\{ a\},\{ b\}) \; :=\; 
\left\{
\begin{array}{ll}
1 & \mb{if $\{ a\} = \{ b\}$,} \\
0 & \mb{if $\{ a\} \neq \{ b\}$.} \\
\end{array}
\right.
\]
Let $C_a\leq\Sl_n$ denote the column stabilizer of $[a]$. A {\it $\lambda$-polytabloid} is given by
\[
\spi{a} \; := \; \sum_{\sigma\in C_a} \{ a\}\sigma\eps_\sigma\;\in\; M^\lambda\; .
\]
The {\it Specht module} $S^\lambda$ is defined to be the $\Z$-linear span of the $\lambda$-polytabloids in $M^\lambda$. It carries a $\Z\Sl_n$-module structure 
as a submodule of $M^\lambda$, and it carries an $\Sl_n$-invariant bilinear form by restriction of $(-,=)$ to $S^\lambda$, again denoted by $(-,=)$. To simplify 
notation, we sometimes rescale by the James factor to
\[
(-,=)_0 \; :=\; \left(\prod_{i\geq 1} (\lambda'_i - \lambda'_{i+1})!\right)^{\!\! -1}\cdot (-,=)\; ,
\]
cf.\ {\bf\cite[\rm 10.4]{J}}.

A {\it standard basis} of $S^\lambda$ over $\Z$ is given by the set 
of standard polytabloids, where a standard polytabloid $\spi{a}$ is attached to a tableau $[a]$ with strictly increasing rows from left to right, and strictly 
increasing columns from top to bottom. For the Garnir relations between polytabloids we refer to {\bf\cite[\rm 7.2]{J}}. Let
\[
\begin{array}{rcl}
S^\lambda & \lraa{\eta} & S^{\lambda,\ast} \\
\xi       & \lramaps    & (\xi,-)\; , \\
\end{array}
\]
whose cokernel shall be denoted by $S^{\lambda,\ast}/S^\lambda$. 

For a commutative ring $A$, we denote by $S_A^\lambda$ the $A\Sl_n$-module $A\ts_\sZ S^\lambda$. 

Given a prime $p$, the Jantzen filtration of $S_{\sF_p}^\lambda$ with respect to $S_{\sZ_{(p)}}^\lambda\lraa{\eta} S_{\sZ_{(p)}}^{\lambda,\ast}$ has been 
defined in (\ref{DefH0}) as being given by
\[
S_{\sF_p}^\lambda(i) \; =\; \Big((p^i S_{\sZ_{(p)}}^{\lambda,\ast})\eta^{-1} + pS_{\sZ_{(p)}}\Big) / pS_{\sZ_{(p)}}\;\tm\; S_{\sF_p}^\lambda\; ,
\]
yielding a decreasing filtration as $i$ runs over $\Z_{\geq 0}$. Note that $S_{\sF_p}^\lambda(0)/S_{\sF_p}^\lambda(1)\iso D_{\sF_p}^\lambda$, which is, 
according to the convention adopted here, nonzero if and only if $\lambda$ is $p$-regular {\bf\cite[\rm 12.2, 11.1]{J}}. In this case, $D_{\sF_p}^\lambda$ does 
not appear in $S_{\sF_p}^\lambda(1)$ by {\bf\cite[\rm 12.2]{J}}.

In a linear combination expressing $[S_{\sF_p}^\lambda]$ in $K_0(\modfr\Z_{(p)}\Sl_n)$ in terms of simple modules, a lower index $i\geq 0$ 
indicates that this summand appears as a summand of the subquotient $[S_{\sF_p}^\lambda(i)/S_{\sF_p}^\lambda(i + 1)]$. Thus, this summand appears with 
multiplicity $i$ in $[S_{\sZ_{(p)}}^{\lambda,\ast}/S_{\sZ_{(p)}}^\lambda]$ (\ref{LemN1} iii). This is similar to the notation in {\bf\cite{L}}, where one finds 
e.g.\
\[
\begin{array}{rcl}
[S_{\sF_2}^{(4,3,2^2,1)}]
& = & [D_{\sF_2}^{(8,3,1)}]_2 + [D_{\sF_2}^{(6,5,1)}]_3 + [D_{\sF_2}^{(10,2)}]_4 + [D_{\sF_2}^{(10,2)}]_7 + 3[D_{\sF_2}^{(12)}]_4 + 2[D_{\sF_2}^{(12)}]_7 \\
&   & + [D_{\sF_2}^{(6,4,2)}]_5 + [D_{\sF_2}^{(6,4,2)}]_8 + [D_{\sF_2}^{(7,3,2)}]_5 + [D_{\sF_2}^{(8,4)}]_5 + [D_{\sF_2}^{(7,5)}]_6 + [D_{\sF_2}^{(7,5)}]_7 \\
&   & + [D_{\sF_2}^{(11,1)}]_6 + 4[D_{\sF_2}^{(11,1)}]_9 + 2[D_{\sF_2}^{(9,3)}]_7 + [D_{\sF_2}^{(9,3)}]_9  + [D_{\sF_2}^{(5,4,2,1)}]_{10}\; . \\
\end{array}
\]
(Presumably, there are also examples in which the multiplicity of a simple module is distributed over more than two Jantzen subquotients.)
Thus for instance,
\[
[S_{\sF_2}^{(4,3,2^2,1)}(7)/S_{\sF_2}^{(4,3,2^2,1)}(8)] \; =\; [D_{\sF_2}^{(10,2)}] + 2[D_{\sF_2}^{(12)}] + [D_{\sF_2}^{(7,5)}] + 2[D_{\sF_2}^{(9,3)}]\; . 
\]
Hence the multiplicity of $\Z/2^7\Z$ as a summand of the abelian group $S_{\sZ_{(2)}}^{(4,3,2^2,1),\ast}/S_{\sZ_{(2)}}^{(4,3,2^2,1)}$ is given by 
$\dim_{\sF_2} D_{\sF_2}^{(10,2)} + 2\dim_{\sF_2} D_{\sF_2}^{(12)} + \dim_{\sF_2} D_{\sF_2}^{(7,5)} + 2\dim_{\sF_2} D_{\sF_2}^{(9,3)}$ (\ref{RemN0}).

\subsection{Two-row partitions}

\begin{quote}
\begin{footnotesize}
We consider Specht modules indexed by two-row partitions, that is, partitions of the form $(n-m,m)$.
The basic ingredients that allow to use Jantzen's Lemma are the theorem of {\sc James} on the decomposition numbers of Specht modules of two-row partitions, 
which are contained in $\{ 0,1\}$, and the theorem of {\sc Schaper,} expressing $[S_{\sZ_{(p)}}^{\lambda,\ast}/S_{\sZ_{(p)}}^\lambda]$ as linear combination of 
Specht modules $[S_{\sF_p}^\lambda]$. 
\end{footnotesize}
\end{quote}

Let $n\geq 1$, let $0\leq m\leq n/2$ and let $p$ be a prime number. 

\begin{Lemma}[{\bf\cite[20.1]{J}}]
\label{LemT1}
We have $\;\rk_\sZ\, S^{(n-m,m)} = \smatze{n}{m}\,\cdot\,\frac{n-2m+1}{n-m+1}$.
\end{Lemma}

Suppose given integers $s$ and $t$. If $s\geq 0$ and $t\geq 1$, we write them $p$-adically as $s+1 = \sum_{i\in [0,k]} s_i p^i$ and 
$t = \sum_{i\in [0,l]} t_i p^i$, where $s_i, t_i\in [0,p-1]$ and where $s_k\neq 0$ and $t_l\neq 0$. The integer $s+1$ is said to {\it contain $t$ to base $p$} 
if ($s\geq 0$, $t\geq 1$, $k > l$ and $t_i\in \{ 0, s_i\}$ for $i\in [0,l]$) or if ($s\geq 0$ and $t = 0$). Let 
\[
f_p(s,t) \; :=\; 
\left\{
\begin{array}{ll}
1 & \mb{if $s+1$ contains $t$ to base $p$} \\
0 & \mb{otherwise}\; . \\
\end{array}
\right.
\]

\begin{Theorem}[{\bf\cite[24.15]{J}}]
\label{ThT2}
Let $j\in [1,m]$. We have 
\[
[S_{\sF_p}^{(n-m,m)} : D_{\sF_p}^{(n-j,j)}] \; =\; f_p(n-2j,m-j)\; (\in \{ 0,1\})\; . 
\]
\end{Theorem}

\begin{Theorem}[particular case of Schaper's formula {\bf\cite[\rm p.\ 60]{Sch}}, see also {\bf\cite[5.33]{Mathas}}]
\label{ThT3}\Absatz
In $K_0(\modfr\Z_{(p)}\Sl_n)$, we have 
\[
[S_{\sZ_{(p)}}^{(n-m,m),\ast}/S_{\sZ_{(p)}}^{(n-m,m)}] \; =\; \sumd{i\in [0,m-1]} v_p\left({\text\frac{n+1-m-i}{m-i}}\right) [S_{\sF_p}^{(n-i,i)}]\; . 
\]
\end{Theorem}

For $j\in [0,m]$, we abbreviate 
\[
\mu(n,m,p\, ;\, j) \; :=\; \sum_{i\in [j,m-1]}v_p\!\left(\frac{n+1-m-i}{m-i}\right) f_p(n-2j,i-j)\; . 
\]

\begin{Corollary}
\label{CorT3_5}
Combining (\ref{ThT2}) and (\ref{ThT3}), we get
\[
[S_{\sZ_{(p)}}^{(n-m,m),\ast}/S_{\sZ_{(p)}}^{(n-m,m)}]\; =\; \sumd{j\in [0,m-1]} \mu(n,m,p\, ;j)\, [D^{(n-j,j)}]\; .
\]
\end{Corollary}

We denote by $(b_{p\, ;\, i,j})_{0\leq i,j\leq n/2}$ the inverse of the lower triangular unipotent integral matrix $(f_p(n-2j,i-j))_{0\leq i,j\leq n/2}\;$, and 
remark that $\dim_{\sF_p} D_{\sF_p}^{(n-j,j)} = \sum_{l\in [0,j]} b_{p\, ;\, j,l} \smatze{n}{l}\cdot\frac{n-2l+1}{n-l+1}$.

\begin{Theorem}
\label{CorT4}
In $K_0(\modfr\Z_{(p)}\Sl_n)$, we have
$$
[S_{\sF_p}^{(n-m,m)}] \; =\; \sumd{j\in [0,m]} f_p(n-2j,m-j) [D_{\sF_p}^{(n-j,j)}]_{\mu(n,m,p\, ;\, j)}\; .
\leqno (\ref{CorT4}.1)
$$
Thus as $\Z_{(p)}$-modules, we have
$$
S_{\sZ_{(p)}}^{(n-m,m),\ast}/S_{\sZ_{(p)}}^{(n-m,m)}\;\iso\; 
{\dis\Ds_{j\,\in\, [0,m-1]}}\left(\Z/p^{\mu(n,m,p\, ;\, j)}\Z\right)^{\sum_{l\in [0,j]} b_{p\, ;\, j,l} \ssmatze{n}{l}\cdot\frac{n-2l+1}{n-l+1}}\; . \\
\leqno (\ref{CorT4}.2)
$$
In other words, the right hand side lists the $p$-part of the elementary divisors of the Gram matrix of the invariant bilinear form on the Specht module 
$S^{(n-m,m)}$, which is unique up to scalar, in an unordered manner. 

\rm
Formula (\ref{CorT4}.1) follows by an application of (\ref{CorN2}) to (\ref{ThT2}) and (\ref{CorT3_5}), where $R = \Z_{(p)}$, $\pi = p$, 
$\Lambda = \Z_{(p)}\Sl_n$, and where $(X\lraa{\phi} Y) = (S_{\sZ_{(p)}}^{(n-m,m)}\lraa{\eta} S_{\sZ_{(p)}}^{(n-m,m),\ast})$. Cf.\ {\bf\cite[\rm 4.12, 11.5]{J}}.

Formula (\ref{CorT4}.2) now ensues by (\ref{RemN0}).
\end{Theorem}

We draw a conclusion that will follow again from (\ref{ThL4}) below.

\begin{Corollary}
\label{CorT5}
If $p > m$, we obtain 
$$
[S_{\sF_p}^{(n-m,m)}] \; =\; [D_{\sF_p}^{(n-m,m)}]_0 + \sumd{j\in [0,m-1]} \{ n\con_p m+j-1\} [D_{\sF_p}^{(n-j,j)}]_{v_p(n+1-m-j)}\; .
\leqno (\ref{CorT5}.1)
$$
Thus as $\Z_{(p)}$-modules, we have
$$
S_{\sZ_{(p)}}^{(n-m,m),\ast}/S_{\sZ_{(p)}}^{(n-m,m)}\;\iso\; 
{\dis\Ds_{j\,\in\, [0,m-1]}} \Big(\Z_{(p)}/(n+1-m-j)\Z_{(p)}\Big)^{\ssmatze{n}{j}\cdot\frac{n-2j+1}{n-j+1}}\; . 
\leqno (\ref{CorT5}.2)
$$

\rm
For $j\in [0,m-1]$, we obtain $\mu_{n,m,p\, ;j} = v_p(n+1-m-j)$, whence (\ref{CorT5}.1).

In (\ref{CorT5}.2) and (\ref{CorT4}.2), we remark that the summands are zero unless $j\con_p n + 1 - m$, which happens at most once. If $j\con_p n + 1 - m$, we 
recall that the outer exponent in (\ref{CorT4}.2) is just $\dim_{\sF_p} D_{\sF_p}^{(n-j,j)}$, which equals $\dim_{\sF_p} S_{\sF_p}^{(n-j,j)}$ by (\ref{ThT2}).
\end{Corollary}

\subsection{Numerical examples}
\label{SubSecNumEx}

\begin{footnotesize}
\begin{quote}
\begin{Example}[case $m = 1$]
\label{ExT6_5}\rm
Let $n\geq 2$. For $p$ arbitrary, (\ref{CorT5}.2) yields
\[
S_{\sZ_{(p)}}^{(n-1,1),\ast}/S_{\sZ_{(p)}}^{(n-1,1)}\;\iso\; \Z_{(p)}/n\Z_{(p)}\; .
\]
This is also a particular case of (\ref{ThHook}) below.
\end{Example}

\begin{Example}[case $m = 2$]
\label{ExT7}\rm
Let $n\geq 4$. For $p > 2$, (\ref{CorT5}.2) yields
\[
S_{\sZ_{(p)}}^{(n-2,2),\ast}/S_{\sZ_{(p)}}^{(n-2,2)}\;\iso\; \Big(\Z_{(p)}/(n-1)\Z_{(p)}\Big) \;\ds\; \Big(\Z_{(p)}/(n-2)\Z_{(p)}\Big)^{n-1}\; .
\]

For $p = 2$, (\ref{CorT4}.2) yields
\[
S_{\sZ_{(2)}}^{(n-2,2),\ast}/S_{\sZ_{(2)}}^{(n-2,2)}\;\iso\; \Big(\Z_{(2)}/\smatze{n-1}{2}\Z_{(2)}\Big) \;\ds\; \Big(\Z_{(2)}/(n-2)\Z_{(2)}\Big)^{n - 2}\; .
\]
Note that for $n = 4$, where the elementary divisors of $S^{(2,2)}$ are $2$ and $6$, the result is formulated using redundant zero summands.
In (\ref{Th2^2}) below, we will give two bases essentially diagonalizing the Gram matrix of the Specht module to the transposed partition (cf.\ \ref{CorTr3}).
\end{Example}

\begin{Example}[case $m = 3$]
\label{ExT8}\rm
Let $n\geq 6$. For $p > 3$, (\ref{CorT5}.2) yields
\[
S_{\sZ_{(p)}}^{(n-3,3),\ast}/S_{\sZ_{(p)}}^{(n-3,3)} \iso 
        \Big(\Z_{(p)}/(n-2)\Z_{(p)}\Big) 
  \ds   \Big(\Z_{(p)}/(n-3)\Z_{(p)}\Big)^{n-1}
  \ds   \Big(\Z_{(p)}/(n-4)\Z_{(p)}\Big)^{\frac{n(n-3)}{2}}\; .
\]

For $p = 3$, (\ref{CorT4}.2) yields
\[
\begin{array}{rcl}
S_{\sZ_{(3)}}^{(n-3,3),\ast}/S_{\sZ_{(3)}}^{(n-3,3)}
& \iso &         \Big(\Z_{(3)}/\smatze{n-2}{3}\Z_{(3)}\Big) 
         \;\ds\; \Big(\Z_{(3)}/(n-3)\Z_{(3)}\Big)^{n - 2} \\
&      & \;\ds\; \Big(\Z_{(3)}/(n-4)\Z_{(3)}\Big)^{\frac{n(n-3)}{2} - 1}\; . \\
\end{array}
\] 

For $p = 2$, (\ref{CorT4}.2) yields
\[
\begin{array}{rcl}
S_{\sZ_{(2)}}^{(n-3,3),\ast}/S_{\sZ_{(2)}}^{(n-3,3)}
& \iso &          \Big(\Z_{(2)}/\frac{n-2}{2^{\{n\con_4 0\}}}\Z_{(2)}\Big) 
         \;\ds\; \Big(\Z_{(2)}/\smatze{n-3}{2}\Z_{(2)}\Big)^{(n - 1) - \{n\con_2 0\}} \\
&      & \;\ds\; \Big(\Z_{(2)}/(n-4)\Z_{(2)}\Big)^{\frac{n(n-3)}{2} - (n - 1) + \{n\con_4 0\}}\; . \\
\end{array}
\] 
\end{Example}

\begin{Example}[case $m = 4$]
\label{ExT9}\rm
Let $n\geq 8$. For $p > 4$, (\ref{CorT5}.2) yields
\[
\begin{array}{rcl}
S_{\sZ_{(p)}}^{(n-4,4),\ast}/S_{\sZ_{(p)}}^{(n-4,4)} 
& \iso &        \Big(\Z_{(p)}/(n-3)\Z_{(p)}\Big) 
         \;\ds\;\Big(\Z_{(p)}/(n-4)\Z_{(p)}\Big)^{n-1} \\
&      & \;\ds\;\Big(\Z_{(p)}/(n-5)\Z_{(p)}\Big)^{\frac{n(n-3)}{2}}
         \;\ds\;\Big(\Z_{(p)}/(n-6)\Z_{(p)}\Big)^{\frac{n(n-1)(n-5)}{6}}\; .
\end{array}
\]

For $p = 3$, (\ref{CorT4}.2) yields
\[
\begin{array}{rcl}
S_{\sZ_{(3)}}^{(n-4,4),\ast}/S_{\sZ_{(3)}}^{(n-4,4)}
& \iso &        \Big(\Z_{(3)}/\left(3^{-\{n\con_9 0,6\}}(n-3)\right)\Z_{(3)}\Big) 
         \;\ds\;\Big(\Z_{(3)}/\smatze{n-4}{3}\Z_{(3)}\Big)^{(n-1) - \{ n\con_3 0\}} \\
&      & \;\ds\;\Big(\Z_{(3)}/(n-5)\Z_{(3)}\Big)^{\frac{n(n-3)}{2} - (n-1)} \\
&      & \;\ds\;\Big(\Z_{(3)}/(n-6)\Z_{(3)}\Big)^{\frac{n(n-1)(n-5)}{6} - (n-1) + \{n\con_9 0,6\}} \; . \\
\end{array}
\] 

For $p = 2$, (\ref{CorT4}.2) yields
\[
\begin{array}{rcl}
S_{\sZ_{(2)}}^{(n-4,4),\ast}/S_{\sZ_{(2)}}^{(n-4,4)}
& \iso &        \Big(\Z_{(2)}/\smatze{n-3}{4}\Z_{(2)}\Big) 
         \;\ds\;\Big(\Z_{(2)}/\left(2^{-\{ n\con_4 2\}}(n-4)\right)\Z_{(2)}\Big)^{n-2} \\
&      & \;\ds\;\Big(\Z_{(2)}/\smatze{n-5}{2}\Z_{(2)}\Big)^{\frac{n(n-3)}{2} - \{ n\con_4 2\}(n-1) - \{ n\con_4 1\}} \\
&      & \;\ds\;\Big(\Z_{(2)}/(n-6)\Z_{(2)}\Big)^{\frac{n(n-1)(n-5)}{6} - \frac{n(n-3)}{2} + \{ n\con_4 2\}(n-2)} \; . \\
\end{array}
\] 
\end{Example}
\end{quote}
\end{footnotesize}

\subsection{Unimodular lattices in $S_\sQ^{(n-m,m)}$}
\label{SubSecUnimod}

\begin{quote}
\begin{footnotesize}
Closely connected to the question for the elementary divisors of the Gram matrix is the question of the existence of unimodular lattices in a Specht module. For 
two-row partitions this investigation has been initiated by {\sc Plesken} {\bf\cite{P93}}. Our result now asserts the non-existence of such lattices in the 
cases not treated in loc.\ cit. Reformulated, this amounts to the assertion that a certain system of two Pell equations is only trivially solvable, the proof of 
which we owe to {\sc E.\ Wirsing.} The number of solutions of general Pellian systems has been studied extensively (cf.\ {\bf\cite{Bennett}}, where also further 
references may be found). The method employed in the particular case here does not seem to generalize.
\end{footnotesize}
\end{quote}

Let $R$ be a localization of $\Z$ at a maximal ideal, or $\Z$ itself, let $\lambda$ be a partition of $n$, and let $X\tm S_\sQ^\lambda$ be a full 
$R\Sl_n$-lattice, i.e.\ an $R\Sl_n$-submodule that is finitely generated free over $R$ of rank $\dim_\sQ S_\sQ^\lambda$. We denote by
\[
X^\# \; :=\; \{ v\in S_\sQ^\lambda\; |\; (v,X)\tm R\}\;\tm\; S_\sQ^\lambda
\]
its dual lattice. Note that $X^\# \iso X^\ast$, the latter denoting the abstract $R$-dual. 

If $X\iso X^\ast$, the $R\Sl_n$-lattice $X$ is called {\it unimodular.} Under our assumptions, this is equivalent to the existence of a scalar $a\in\Q$ such 
that $X = a X^\#$. 

Given an inclusion $X\tm Y\tm S_\sQ^\lambda$ of full $R\Sl_n$-lattices, by self duality of simple $\F_p\Sl_n$-modules for each prime $p$, we have 
$[Y/X] = [X^\#/Y^\#]$ in the 
Grothendieck group {\bf\cite[\rm 11.5]{J}}.

We include a proof of a corollary of {\sc Plesken,} restricting his argument to this corollary as well.

\begin{Proposition}[{\bf\cite[Cor.\ II.4]{P93}}]
\label{PropU1}
Let $p$ be a prime and assume $\lambda$ to be $p$-regular. If $S_\sQ^\lambda$ contains a unimodular $\Z_{(p)}\Sl_n$-lattice, then 
\[
[S_{\sZ_{(p)}}^{\lambda,\ast}/S_{\sZ_{(p)}}^\lambda]\;\in\; 2 K_0(\modfr\Z_{(p)}\Sl_n)\; . 
\]

\rm
Let $X\tm S_\sQ^\lambda$ be a unimodular $\Z_{(p)}\Sl_n$-lattice, let $a\in\Q$ such that $X = a X^\#$. Let $m\in\Z_{(p)}$ be such that 
$m S_{\sZ_{(p)}}^\lambda\tm X$. The filtration
\[
m S_{\sZ_{(p)}}^\lambda\;\tm\; X = a X^\#\;\tm\; am^{-1} S_{\sZ_{(p)}}^{\lambda,\#}
\]
shows that $[am^{-1} S_{\sZ_{(p)}}^{\lambda,\#}/m S_{\sZ_{(p)}}^\lambda]\in 2 K_0(\modfr\Z_{(p)}\Sl_n)$. Therefore, there exists $s\in\Z$ such that
\[
[S_{\sZ_{(p)}}^{\lambda,\ast}/S_{\sZ_{(p)}}^\lambda] + s [S_{\sF_p}^\lambda]\;\in\; 2 K_0(\modfr\Z_{(p)}\Sl_n)
\] 
Since $\lambda$ is $p$-regular, counting multiplicities of $D_{\sF_p}^\lambda$ yields $s\in 2\Z$ by (\ref{LemN1} i) and {\bf\cite[\rm 12.2]{J}}.
\end{Proposition}  

\begin{Remark}
\label{RemU2}
The converse to (\ref{PropU1}) holds as well, provided the decomposition numbers of $S_{\sF_p}^\lambda$ are in $\{ 0,1\}$. For this direction, $\lambda$ need 
not be $p$-regular.

\rm
In fact, suppose $[S_{\sZ_{(p)}}^{\lambda,\ast}/S_{\sZ_{(p)}}^\lambda]\;\in\; 2 K_0(\modfr\Z_{(p)}\Sl_n)$. Let 
$S_{\sZ_{(p)}}^\lambda\tm M\tm S_\sQ^{\lambda,\ast}$ be a full $\Z_{(p)}\Sl_n$-lattice that is maximal with respect to the property that $(M,M)\tm \Z_{(p)}$. 
Let $a\geq 0$ such that $p^{2^a} M^\#\tm M$. 

If $p^{2^b} M^\#\tm M$ for some $b\geq 1$, then $(p^{2^{b-1}} M^\#, p^{2^{b-1}} M^\#) = (M^\#, p^{2^b} M^\#)\tm (M^\#, M)\tm\Z_{(p)}$, whence 
$M = M + p^{2^{b-1}} M^\#$ by maximality of $M$, i.e.\ $p^{2^{b-1}} M^\#\tm M$. By induction, starting with $b = a$, we conclude that $p M^\#\tm M\tm M^\#$. 

Now by the filtration $S_{\sZ_{(p)}}^\lambda\tm M\tm M^\#\tm S_{\sZ_{(p)}}^{\lambda,\#}$, the decomposition numbers of $M^\#/M$ are even. But since $M^\#/M$ is a 
quotient of $M^\#/pM^\#$, the decomposition numbers of which are in $\{ 0,1\}$ by assumption, we infer that $M = M^\#$.
\end{Remark}

Now let $\lambda = (n-m,m)$. The cases $m = 1$ and $m = 2$ have been treated in {\bf\cite[\rm p.\ 98 and II.5]{P93}}. 

\begin{Theorem}
\label{ThU3}
Let $3\leq m\leq n/2$. The module $S_{\sQ}^{(n-m,m)}$ does not contain a unimodular $\Z\Sl_n$-lattice.

\rm
Given a prime $p$, Schaper's formula reads
\[
[S_{\sZ_{(p)}}^{(n-m,m),\ast}/S_{\sZ_{(p)}}^{(n-m,m)}]\; =\; \sumd{i\in [1,m]} v_p\!\left(\frac{n - 2m + 1 + i}{i}\right) [S_{\sF_p}^{(n-m+i,m-i)}] \; .
\]
First, let us consider the case $2m = n$. If there was a unimodular $\Z\Sl_n$-lattice, there would be a unimodular $\Z_{(3)}\Sl_n$-lattice. However, 
$D_{\sF_3}^{(m+2,m-2)}$ appears in $[S_{\sZ_{(3)}}^{(n-m,m),\ast}/S_{\sZ_{(3)}}^{(n-m,m)}]$ with multiplicity $1$, so this is impossible by (\ref{PropU1}).

Now, let us consider the case $2m < n$ and assume the existence of a unimodular $\Z\Sl_n$-lattice. At a prime $p$, existence of a unimodular 
$\Z_{(p)}\Sl_n$-lattice implies by (\ref{PropU1}) that all multiplicities of simple modules in $[S_{\sZ_{(p)}}^{(n-m,m),\ast}/S_{\sZ_{(p)}}^{(n-m,m)}]$ are 
even. Considering the multiplicity of $D_{\sF_p}^{(n-m+i,m-i)}$ for $i\in [1,m]$, beginning with $i = 1$, we obtain the condition that $2$ 
divides $v_p\!\left(\frac{n - 2m + 1 + i}{i}\right)$ for all $i\in [1,m]$. Since this holds for all primes $p$, and at least for $i\in [1,3]$, we conclude that 
there exist positive integers $x,y,z$ such that 
\[
\begin{array}{rcl}
n - 2m + 2 & = & x^2 \\
n - 2m + 3 & = & 2y^2 \\
n - 2m + 4 & = & 3z^2\; , \\
\end{array}
\]
whence
$$
\begin{array}{rcl}
2y^2 - x^2 & = & 1 \\
3z^2 - x^2 & = & 2\; . \\
\end{array}
\leqno (\ast)
$$
{\sc E.\ Wirsing} {\bf\cite{W}} proved that $(x,y,z) = (1,1,1)$ is the only solution of $(\ast)$ in positive integers, as reproduced below. Since $x = 1$ would 
correspond to $n = 2m - 1$, this assertion contradicts the assumption on the existence of a unimodular $\Z\Sl_n$-lattice and proves the theorem.

First, we remark that a solution $(x,y,z)$ consists of pairwise coprime integers, and that $x\con_2 1$.

The rational points $(\xi,\zeta)$ on the ellipse $\xi^2 + 3\zeta^2 = 4$ are parametrized by 
\[
(\xi,\zeta)\; =\; \left(2\frac{1-3t^2}{1+3t^2},\frac{4t}{1+3t^2}\right)\; ,
\]
where $t\in\Q\cup\{\infty\}$. Letting $(\xi,\zeta) = (\frac{x}{y},\frac{z}{y})$, a solution $(x,y,z)$ yields such a rational point. We may exclude $t = \infty$, 
corresponding to $(-2,0)$, since $z\neq 0$. Writing $t = \frac{r}{s}$ with coprime positive integers $r$ and $s$, we 
obtain $r\con_2 s\con_2 1$ since $r\not\con_2 s$ would imply $x\con_2 0$. If $s\not\con_3 0$, then $rs$ and $(s^2 + 3r^2)/4$ are coprime, whence 
\[
(x,y,z)\; =\; \left((s^2 - 3r^2)/2,(s^2 + 3r^2)/4,rs\right)\; .
\]
Now $x^2 = 3z^2 - 2$ yields 
$$
s^4 - 1\; =\; 2\cdot\left(\frac{3}{4}(r^2 - s^2)\right)^2\; .
\leqno (\ast\ast_1)
$$

If $s\con_3 0$, then $rs/3$ and $(s^2 + 3r^2)/12$ are coprime, whence 
\[
(x,y,z)\; =\; \frac{1}{3}\left((s^2 - 3r^2)/2,(s^2 + 3r^2)/4,rs\right)\; .
\]
Now $x^2 = 3z^2 - 2$ yields 
$$
r^4 - 1\; =\; 2\cdot\left(\frac{1}{12}(s^2 - 9r^2)\right)^2\; .
\leqno (\ast\ast_2)
$$
By a result of {\sc Euler,} however, the integral equation
$$
u^4 - v^4\; =\; 2w^2
\leqno (\ast\ast)
$$
is unsolvable if $w\neq 0$ {\bf\cite[\rm p.\ 82]{We}}. Therefore, the only solution to $(\ast)$ in positive integers results from $(\ast\ast_2)$ as being 
$(x,y,z) = (1,1,1)$.
\end{Theorem}
\subsection{Some three- and four-part partitions}
\label{Sec34part}

\begin{quote}
\begin{footnotesize}
The partial results and conjectures that follow we could lift from the table of {\sc F.\ L\"ubeck} {\bf\cite{L}} and from calculations of 
decomposition numbers due to {\sc James}, {\sc Williams,} {\sc To Law,} {\sc Benson,} {\sc M\"uller} et al. 
{\bf\cite{B,J,JW,Atlas,M,M18,TL}}. See also {\bf\cite[\rm p.\ 113]{J}}. A general result for $p$ large that covers the respective first cases listed in this 
section, is given below (\ref{ThL4}).
\end{footnotesize}
\end{quote}

Recall that we stipulated the module $D_{\sF_p}^\mu$ to be zero whenever $\mu$ is not a $p$-regular partition.

\begin{Proposition}[$\mathbf{n-3,2,1}$]
\label{PropS1}
Let $n\geq 5$. For $p > 3$, we have 
\[
[S_{\sF_p}^{(n-3,2,1)}] \; =\; [D_{\sF_p}^{(n-3,2,1)}]_0 + \{n\con_p 3\}[D_{\sF_p}^{(n-2,1^2)}]_{v_p(n-3)} + \{n\con_p 1\}[D_{\sF_p}^{(n-2,2)}]_{v_p(n-1)}\; .
\]
For $p = 2$, we have
\[
\begin{array}{rcl}
[S_{\sF_2}^{(n-3,2,1)}]
& = & [D_{\sF_2}^{(n-3,2,1)}]_0 + \{ n\con_2 1\}[D_{\sF_2}^{(n-2,2)}]_{v_2(n-1) + v_2(n-3)} \\ 
&   & + \{ n\con_2 1\} [D_{\sF_2}^{(n)}]_{2\{n\con_4 1\} + (v_2(n-3) - 1)\{ n\con_4 3\}}\; .\\
\end{array}
\]
For $p = 3$, we have
\[
\begin{array}{rcl}
[S_{\sF_3}^{(n-3,2,1)}] 
& = & [D_{\sF_3}^{(n-3,2,1)}]_0 + [D_{\sF_3}^{(n-3,3)}]_1 + \{ n\con_3 1\}[D_{\sF_3}^{(n-2,2)}]_{1 + v_3(n-1)}\\
&   & + \{ n\con_3 0\}[D_{\sF_3}^{(n-1,1)}]_{1 + v_3(n-3)} + \{ n\con_3 0\}[D_{\sF_3}^{(n-2,1^2)}]_{v_3(n-3)} \\
&   & + (1 + \{n\con_9 2,3,4\})[D_{\sF_3}^{(n)}]_1\; . \\
\end{array} 
\]

\rm
Schaper's formula (see e.g.\ {\bf\cite[\rm 5.33]{Mathas}}) reads
\[
\begin{array}{rcl}
[S_{\sZ_{(p)}}^{(n-3,2,1),\ast}/S_{\sZ_{(p)}}^{(n-3,2,1)}] 
& = & v_p(3) [S_{\sF_p}^{(n-3,3)}] + v_p(n-1) [S_{\sF_p}^{(n-2,2)}] - v_p\!\left(\frac{n-1}{3}\right)[S_{\sF_p}^{(n)}] \\
&   & + v_p(n-3)[S_{\sF_p}^{(n-2,1^2)}]\; ,
\end{array}
\]
where $p$ is an arbitrary prime, and where $n\geq 6$. We expand this formula into simple modules using {\bf\cite[\rm 24.15]{J}} and {\bf\cite[\rm App.]{JW}}. 
Since by loc.\ cit.\ the decomposition numbers of $S_{\sF_p}^{(n-3,2,1)}$ itself are contained in $\{ 0,1,2\}$, we may apply (\ref{CorN2}, \ref{CorN3}) to 
obtain the distribution of the occurring simple
modules over the Jantzen subquotients as stated above. 
\end{Proposition}

\begin{Lemma}
\label{LemS1_5}
Let $n\geq 6$. We have
$\;
[S_{\sF_2}^{(n-4,2^2)} : D_{\sF_2}^{(n-2,2)}]\; =\; \{ n\con_2 0\}\; .
$

\rm
This follows by {\bf\cite[\rm 3.13]{JW}} if $n$ is odd, and by {\bf\cite[\rm 4.5]{JW}} if $n$ is even. We state this explicitly since in this case we do not have 
$\lambda_3 < p$, cf.\ {\bf\cite[\rm Introduction; more precisely, 4.13]{JW}}. 
\end{Lemma}

\begin{Conjecture}
\label{ConjS2}
Let $n\geq 6$. We have
$\;
[S_{\sF_2}^{(n-4,2^2)} : D_{\sF_2}^{(n)}]\; =\; 1 + \{ n\con_4 0,1\}\; .
$
\end{Conjecture}

\begin{Proposition}[$\mathbf{n-4,2^2}$]
\label{PropS3}
Let $n\geq 6$. For $p > 3$, we have
\[
[S_{\sF_p}^{(n-4,2^2)}] \; =\; [D_{\sF_p}^{(n-4,2^2)}]_0 + \{n\con_p 2\}[D_{\sF_p}^{(n-2,2)}]_{v_p(n-2)} +  \{n\con_p 3\}[D_{\sF_p}^{(n-3,2,1)}]_{v_p(n-3)}\; .
\]
Suppose $p = 2$. If (\ref{ConjS2}) holds true, then 
\[
\begin{array}{rcl}
[S_{\sF_2}^{(n-4,2^2)}] 
& = & [D_{\sF_2}^{(n-4,3,1)}]_1 + (1 + \{ n\con_4 0,1\}) [D_{\sF_2}^{(n)}]_{1 + \{n\con_4 2,3\}} \\
&   & + \{ n\con_2 0\} [D_{\sF_2}^{(n-2,2)}]_{v_2(n-2) + 1} + \{ n\con_2 1\} [D_{\sF_2}^{(n-3,2,1)}]_{v_2(n-3) + 1}\; . \\
\end{array}
\]
Suppose $p = 3$. For $n = 6$, we have
\[
[S_{\sF_3}^{(2^3)}] \; =\; [D_{\sF_3}^{(3,2,1)}]_1 + [D_{\sF_3}^{(6)}]_2\; . 
\]
If $n\geq 7$, there exist $a_n, b_n\geq 1$ with $a_n + b_n = (1 + v_3(n-3)) + (1 + v_3(n-3))$ such that 
\[
\begin{array}{rcl}
[S_{\sF_3}^{(n-4,2^2)}] 
& = & [D_{\sF_3}^{(n-4,2^2)}]_0 + [D_{\sF_3}^{(n-4,4)}]_1 + \{n\con_3 0\} [D_{\sF_3}^{(n-3,2,1)}]_{v_3(n-3)} \\ 
&   & + \{n\con_3 0\} [D_{\sF_3}^{(n-3,3)}]_{v_3(n-3)+1} + \{n\con_3 2\} [D_{\sF_3}^{(n-2,2)}]_{v_3(n-2)+1} \\
&   & + (1 + \{ n\con_9 4,5,6 \}) [D_{\sF_3}^{(n-1,1)}]_1 + \{ n\con_9 0,6\}[D_{\sF_3}^{(n)}]_2 \\
&   & + \{n\con_9 3\} ([D_{\sF_3}^{(n)}]_{a_n} + [D_{\sF_3}^{(n)}]_{b_n}) \; . \\
\end{array}
\]

\rm
Schaper's formula reads
\[
\begin{array}{rcl}
[S_{\sZ_{(p)}}^{(n-4,2^2),\ast}/S_{\sZ_{(p)}}^{(n-4,2^2)}] 
& = & v_p\!\left(\frac{3}{2}\right) [S_{\sF_p}^{(n-4,4)}] + v_p\!\left(\frac{n-2}{2}\right) [S_{\sF_p}^{(n-2,2)}]  
      - v_p\!\left(\frac{n-2}{3}\right) [S_{\sF_p}^{(n-1,1)}] \\
&   & + v_p(2) [S_{\sF_p}^{(n-4,3,1)}] + v_p(n-3) [S_{\sF_p}^{(n-3,2,1)}] - v_p\!\left(\frac{n-3}{2}\right) [S_{\sF_p}^{(n-2,1^2)}]\; ,
\end{array}
\]
where $p$ is an arbitrary prime and where $n\geq 8$. The result now follows by {\bf\cite[\rm 24.15]{J}} and {\bf\cite{JW}}, using 
(\ref{CorN2}, \ref{CorN3}, \ref{LemS1_5}). 
\end{Proposition}

\begin{Proposition}[$\mathbf{n-4,3,1}$]
\label{PropS4}
Suppose $n\geq 7$. For $p > 4$, we have
\[
\begin{array}{rcl}
[S_{\sF_p}^{(n-4,3,1)}] 
& = & [D_{\sF_p}^{(n-4,3,1)}]_0 + \{n\con_p 2\}[D_{\sF_p}^{(n-3,3)}]_{v_p(n-2)} + \{n\con_p 4\}[D_{\sF_p}^{(n-2,1^2)}]_{v_p(n-4)} \\
&   & + \{n\con_p 5\}[D_{\sF_p}^{(n-3,2,1)}]_{v_p(n-5)}\; . \\
\end{array}
\]
Suppose $p = 2$. There exist $a_n, b_n, c_n, d_n, e_n, f_n, g_n, h_n, i_n, j_n\geq 1$ with 
\[
\begin{array}{rcl}
a_n + b_n       & = & 1 + (1 + v_2(n-2)) \\ 
c_n + d_n       & = & 2 + 2 \\
e_n + f_n + g_n & = & 1 + 1 + v_2(n-4) \\
h_n + i_n + j_n & = & 1 + 1 + v_2(n-5) \\ 
\end{array}
\]
such that
\[
\begin{array}{rcl}
[S_{\sF_2}^{(n-4,3,1)}] 
& = & [D_{\sF_2}^{(n-4,3,1)}]_0 + [D_{\sF_2}^{(n-4,4)}]_2 + \{n\con_2 0\} [D_{\sF_2}^{(n-3,3)}]_{2 + v_2(n-2)} \\
&   & + \{n\con_2 1\} [D_{\sF_2}^{(n-3,2,1)}]_{v_2(n-5)} + \{n\con_4 0\} [D_{\sF_2}^{(n-2,2)}]_{v_2(n-4)} \\
&   & + \{n\con_4 1\} [D_{\sF_2}^{(n-2,2)}]_{1 + v_2(n-5)} + \{n\con_4 2\} ([D_{\sF_2}^{(n-2,2)}]_{a_n} + [D_{\sF_2}^{(n-2,2)}]_{b_n}) \\
&   & + \{n\con_4 0\} [D_{\sF_2}^{(n-1,1)}]_{2 + v_2(n-4)} + \{n\con_8 2,7\}  [D_{\sF_2}^{(n)}]_2 \\
&   & + \{n\con_8 3,6\} ([D_{\sF_2}^{(n)}]_{c_n} + [D_{\sF_2}^{(n)}]_{d_n}) + 2\{n\con_8 0,1\} [D_{\sF_2}^{(n)}]_1 \\
&   & + \{n\con_8 4\} ([D_{\sF_2}^{(n)}]_{e_n} + [D_{\sF_2}^{(n)}]_{f_n} + [D_{\sF_2}^{(n)}]_{g_n}) \\
&   & + \{n\con_8 5\} ([D_{\sF_2}^{(n)}]_{h_n} + [D_{\sF_2}^{(n)}]_{i_n} + [D_{\sF_2}^{(n)}]_{j_n})\; . \\ 
\end{array}
\]

Suppose $p = 3$. We have 
\[
\begin{array}{rcl}
[S_{\sF_3}^{(n-4,3,1)}] 
& = & [D_{\sF_3}^{(n-4,3,1)}]_0 + \{ n\con_3 2\}[D_{\sF_3}^{(n-3,3)}]_{v_3(n-2) + v_3(n-5)} + \{ n\con_3 2\} [D_{\sF_3}^{(n-3,2,1)}]_{v_3(n-5)} \\
&   & + \{ n\con_3 1\}[D_{\sF_3}^{(n-2,1^2)}]_{v_3(n-4)} + \{ n\con_9 2\}[D_{\sF_3}^{(n)}]_2 + \{ n\con_9 5\}[D_{\sF_3}^{(n)}]_{v_3(n-5) - 1}\; . \\ 
\end{array}
\]

\rm
Schaper's formula reads
\[
\begin{array}{rcl}
[S_{\sZ_{(p)}}^{(n-4,3,1),\ast}/S_{\sZ_{(p)}}^{(n-4,3,1)}] 
& = & - v_p\!\left(\frac{n-2}{4}\right)[S_{\sF_p}^{(n)}] + v_p(n-2)[S_{\sF_p}^{(n-3,3)}] + v_p(4)[S_{\sF_p}^{(n-4,4)}] \\
&   & + v_p\!\left(\frac{n-4}{2}\right)[S_{\sF_p}^{(n-2,1^2)}] + v_p(n-5)[S_{\sF_p}^{(n-3,2,1)}]\; , \\  
\end{array}
\]
where $p$ is an arbitrary prime, and where $n\geq 8$. The result follows by {\bf\cite[\rm 24.15]{J}} and {\bf\cite[\rm App.]{JW}}, using (\ref{CorN2}). 
\end{Proposition}

\begin{Lemma}
\label{LemS5}
Suppose given $n\geq 6$.
\begin{itemize}
\item[(i)] If $n\geq 8$, then $[S_{\sF_2}^{(n-4,2,1^2)} : D_{\sF_2}^{(n-4,3,1)}] = 1$.
\item[(ii)] For $p$ an arbitrary prime, we have $[S_{\sF_p}^{(n-4,2,1^2)} : D_{\sF_p}^{(n-3,2,1)}] = \{n\con_p 1\}$.
\item[(iii)] For $p > 4$ prime, we have $[S_{\sF_p}^{(n-4,2,1^2)} : D_{\sF_p}^{(n-3,1^3)}] = \{n\con_p 4\}$.
\end{itemize}

\rm
This follows from {\bf\cite[\rm 3.6, 3.13]{JW}}. Again, we state this explicitly, this time since \linebreak[4] $(n-4,2,1^2)$ is a four-part-partition. 
Assertion (i) also follows by {\bf\cite[\rm Th.\ A]{JII}}.
\end{Lemma}

\begin{Conjecture}
\label{ConjS6}
Suppose given $n\geq 6$.
\begin{itemize}
\item[(i)] If $n\geq 9$, then $[S_{\sF_2}^{(n-4,2,1^2)} : D_{\sF_2}^{(n-4,4)}] = 1$.
\item[(ii)] If $n\geq 7$, then $[S_{\sF_2}^{(n-4,2,1^2)} : D_{\sF_2}^{(n-3,3)}] = \{ n\con_2 0\}$.
\item[(iii)] We have $[S_{\sF_2}^{(n-4,2,1^2)} : D_{\sF_2}^{(n-2,2)}] = 2 \{ n\con_4 0,1\} + 3\{ n\con_4 2\} + \{ n\con_4 3\}$.
\item[(iv)] We have $[S_{\sF_2}^{(n-4,2,1^2)} : D_{\sF_2}^{(n-1,1)}] = 2 \{ n\con_4 0\} + \{ n\con_4 2\}$.
\item[(v)] If $n\geq 7$, then $[S_{\sF_2}^{(n-4,2,1^2)} : D_{\sF_2}^{(n)}] = \{ n\con_8 7\} + 2\{ n\con_8 0,2,3\} + 3\{n\con_8 1,4,6\} + 4\{ n\con_8 5\}$.
\item[(vi)] We have $[S_{\sF_3}^{(n-4,2,1^2)} : D_{\sF_3}^{(n-3,3)}] = \{ n\con_3 1\}$.
\item[(vii)] We have $[S_{\sF_3}^{(n-4,2,1^2)} : D_{\sF_3}^{(n)}] = \{ n\con_3 1\} + \{ n\con_9 4\}$.
\end{itemize}
\end{Conjecture}

\begin{Proposition}[$\mathbf{n-4,2,1^2}$]
\label{PropS7}
Suppose $n\geq 6$. For $p > 4$, we have
\[
[S_{\sF_p}^{(n-4,2,1^2)}] \; =\; [D_{\sF_p}^{(n-4,2,1^2)}]_0 + \{ n\con_p 1\} [D_{\sF_p}^{(n-3,2,1)}]_{v_p(n-1)} 
+ \{ n\con_p 4\} [D_{\sF_p}^{(n-3,1^3)}]_{v_p(n-4)} \; .
\]

Suppose $p = 2$. For $n = 6$, we obtain
\[
[S_{\sF_2}^{(2^2,1^2)}]\; =\; [D_{\sF_2}^{(4,2)}]_3 + [D_{\sF_2}^{(5,1)}]_2 + [D_{\sF_2}^{(6)}]_2\; .
\]
If $n\geq 7$ and if (\ref{ConjS6} i-v) and (\ref{ConjS2}) hold true, then there exist $a_n, \dots, z_n\geq 1$ with 
\[
\begin{array}{rcl}
a_n + b_n             & = & 1 + (3 + v_2(n-4)) \\ 
c_n + d_n + e_n       & = & 1 + 1 + 4 \\
f_n + g_n             & = & 3 + 3 \\
h_n + i_n + j_n       & = & 2 + 2 + (1 + v_2(n-1)) \\
k_n + l_n             & = & 2 + 2 \\
m_n + o_n + q_n       & = & 2 + 2 + (1 + v_2(n-4)) \\
r_n + s_n + t_n + u_n & = & 1 + 1 + 3 + 3 \\
v_n + w_n + x_n       & = & 1 + 1 + 3 \\
y_n + z_n             & = & (1 + v_2(n-4)) + (1 + v_2(n-4)) \\
\end{array}
\]
such that
\[
\begin{array}{rcl}
[S_{\sF_2}^{(n-4,2,1^2)}] 
& = & [D_{\sF_2}^{(n-4,3,1)}]_3 + [D_{\sF_2}^{(n-4,4)}]_1 + \{ n\con_2 1\} [D_{\sF_2}^{(n-3,2,1)}]_{3 + v_2(n-1)} \\
&   & + \{ n\con_2 0\} [D_{\sF_2}^{(n-3,3)}]_{1 + v_2(n-4)} + (2\{ n\con_4 1\} + \{ n\con_4 3\})[D_{\sF_2}^{(n-2,2)}]_1 \\
&   & + \{ n\con_4 0\}([D_{\sF_2}^{(n-2,2)}]_{a_n} + [D_{\sF_2}^{(n-2,2)}]_{b_n}) \\
&   & + \{ n\con_4 2\}([D_{\sF_2}^{(n-2,2)}]_{c_n} + [D_{\sF_2}^{(n-2,2)}]_{d_n} + [D_{\sF_2}^{(n-2,2)}]_{e_n}) \\
&   & + \{n\con_8 0\}([D_{\sF_2}^{(n)}]_{f_n} + [D_{\sF_2}^{(n)}]_{g_n}) \\
&   & + \{n\con_8 1\}([D_{\sF_2}^{(n)}]_{h_n} + [D_{\sF_2}^{(n)}]_{i_n} + [D_{\sF_2}^{(n)}]_{j_n}) \\
&   & + \{n\con_8 2,3\}([D_{\sF_2}^{(n)}]_{k_n} + [D_{\sF_2}^{(n)}]_{l_n}) \\
&   & + \{n\con_8 4\}([D_{\sF_2}^{(n)}]_{m_n} + [D_{\sF_2}^{(n)}]_{o_n} + [D_{\sF_2}^{(n)}]_{q_n}) \\
&   & + \{n\con_8 5\}([D_{\sF_2}^{(n)}]_{r_n} + [D_{\sF_2}^{(n)}]_{s_n} + [D_{\sF_2}^{(n)}]_{t_n} + [D_{\sF_2}^{(n)}]_{u_n}) \\
&   & + \{n\con_8 6\}([D_{\sF_2}^{(n)}]_{v_n} + [D_{\sF_2}^{(n)}]_{w_n} + [D_{\sF_2}^{(n)}]_{x_n}) + \{n\con_8 7\}[D_{\sF_2}^{(n)}]_3  \\
&   & + \{n\con_4 2\} [D_{\sF_2}^{(n-1,1)}]_2 + \{n\con_4 0\} ([D_{\sF_2}^{(n-1,1)}]_{y_n} + [D_{\sF_2}^{(n-1,1)}]_{z_n})\; .   \\
\end{array}
\]

Suppose $p = 3$. If (\ref{ConjS6} vi, vii) hold true, then
\[
\begin{array}{rcl}
[S_{\sF_3}^{(n-4,2,1^2)}] 
& = & [D_{\sF_3}^{(n-4,2,1^2)}]_0 + \{ n\con_3 1\} [D_{\sF_3}^{(n-3,3)}]_{v_3(n-1)} + (\{ n\con_3 1\} + \{ n\con_9 4\})[D_{\sF_3}^{(n)}]_{v_3(n-1)} \\
&   & + \{ n\con_3 1\} [D_{\sF_3}^{(n-3,2,1)}]_{v_3(n-1) + v_3(n-4)}\; . \\
\end{array}
\]

\rm
Schaper's formula reads
\[
\begin{array}{rcl}
[S_{\sZ_{(p)}}^{(n-4,2,1^2),\ast}/S_{\sZ_{(p)}}^{(n-4,2,1^2)}] 
& = & v_p(2) [S_{\sF_p}^{(n-4,2^2)}] + v_p(4) [S_{\sF_p}^{(n-4,3,1)}] + v_p(n-1)[S_{\sF_p}^{(n-3,2,1)}] \\
&   & - v_p(2)[S_{\sF_p}^{(n-4,4)}] - v_p\!\left(\frac{n-1}{2}\right)[S_{\sF_p}^{(n-2,2)}] + v_p\!\left(\frac{n-1}{4}\right)[S_{\sF_p}^{(n)}]\\
&   & + v_p(n-4)[S_{\sF_p}^{(n-3,1^3)}] \; ,
\end{array}
\]
where $p$ is an arbitrary prime, and where $n\geq 8$. The result follows by {\bf\cite[\rm 24.15]{J}} and {\bf\cite[\rm App.]{JW}}, using 
(\ref{CorN2}, \ref{CorN3}). 
\end{Proposition}

As a general pattern, one seems to observe the following. 

\begin{Conjecture}
\label{ConjS6_5}
Let $p$ be a prime, let $n\geq 1$. Given a partition $\lambda$ of $n$ and an integer $k\geq 0$, we denote by $\lambda+k$ the partition 
$(\lambda_1 + k,\lambda_2,\lambda_3,\dots)$. Suppose $\lambda$ to be a partition of $n$, and $\mu$ to be a $p$-regular partition of $n$. 
\begin{itemize}
\item[(i)] There exists integers $a,J\geq 0$ such that for $j\geq J$, we have
\[
[S_{\sF_p}^\lambda : D_{\sF_p}^\mu]\; =\; [S_{\sF_p}^{\lambda + jp^a} : D_{\sF_p}^{\mu + jp^a}]\; .
\]
\item[(ii)] Let $\lambda^\mb{\scr\rm r}$ denote the $p$-regularization of $\lambda$ in the sense of {\bf\cite{JII}}, and suppose 
$(\lambda^\mb{\scr\rm r})_1 = \lambda_1$. The multiplicity 
$[S_{\sZ_{(p)}}^{\lambda+k,\ast}/S_{\sZ_{(p)}}^{\lambda+k} : D_{\sF_p}^{(\lambda + k)^\mb{\sscr\rm r}}]$ does not depend on $k\geq 0$. 
(Cf.\ {\bf\cite[\rm Th.\ A]{JII}} and (\ref{CorN2}).)
\end{itemize}
\end{Conjecture}
\section{At a large prime}

\begin{Notation}
\label{NotL1_5}\rm
Let $n\geq 1$. Given a partition $\mu$ of $n$ and $j\in [1,\mu_2]$, we denote by $\mu[j]$ the partition of $n$ defined by
\[
\mu[j]_i \; : =\;
\left\{
\begin{array}{ll}
\mu_1 + \mu_{\mu'_j} - j + 1 & \mb{if $i = 1$} \\
j - 1                        & \mb{if $i = \mu'_j$} \\
\mu_i                        & \mb{otherwise} \\
\end{array}
\right.
\]
for $i\geq 0$. I.e.\ we cut the last row that meets the $j$th column in the $j$th column and append that piece in the first row.

Let $R(\mu) = \{ i\geq 1\; |\; \mu_i > \mu_{i+1}\}$ be the set of row numbers of removable nodes. Given $k\in R(\mu)$, we denote by $\mu(k)$ the partition of 
$n - 1$ defined by
\[
\mu(k)_i \; : =\;
\left\{
\begin{array}{ll}
\mu_i - 1 & \mb{if $i = k$} \\
\mu_i     & \mb{else} \\
\end{array}
\right.
\]
for $i\geq 0$. I.e.\ the diagram of $\mu(k)$ is obtained from the diagram of $\mu$ by removing the last node in row $k$.

Let $\lambda$ be a partition of $n$, and suppose $\lambda\neq (n)$. Suppose $p$ to be a prime strictly bigger than $n - \lambda_1 = \sum_{i\geq 2} \lambda_i$.
Given $j\in [1,\lambda_2]$, we denote by $h_j := \lambda_1 - j + \lambda'_j$ the hook length of the node $(1,j)$ of the diagram of $\lambda$. Note that $p$ 
divides at most one of the numbers $h_j$ for $j\in [1,\lambda_2]$ since $h_1 - h_{\lambda_2} < p$. If $p$ divides $h_t$, we denote $s = \lambda'_t$. Note that 
$t\in [\lambda_{s+1}+1,\lambda_s]$.
\end{Notation}

We shall need a particular case of the direction of the Carter Conjecture that has been proven by {\sc James} and {\sc Murphy.}

\begin{Theorem}[{\bf\cite[p.\ 222]{JM}}]
\label{ThL1}
If $p$ does not divide $h_j$ for all $j\in [1,\lambda_2]$, then $S_{\sF_p}^\lambda$ is irreducible.
\end{Theorem}

The following is a particular case of a result of {\sc Carter} and {\sc Payne.} 

\begin{Theorem}[{\bf\cite[p.\ 425]{CP}}]
\label{ThL2}
Suppose that $p$ divides $h_t$. Then $\Hom_{\sF_p\Sl_n}(S_{\sF_p}^{\lambda[t]},S_{\sF_p}^\lambda) \neq 0$. In particular, 
$[S_{\sF_p}^\lambda : D_{\sF_p}^{\lambda[t]}] \geq 1$.
\end{Theorem}

The idea of proof and the main ingredient of (\ref{PropL3}), hence of (\ref{ThL4}), are due to {\sc Klesh\-chev.}

\begin{Proposition}
\label{PropL3}
Suppose that $p$ divides $h_t$. Then $[S_{\sF_p}^\lambda] =  [D_{\sF_p}^\lambda] + [D_{\sF_p}^{\lambda[t]}]$.

\rm
{\it Consider the case $\lambda_1 > \lambda_2$.}
By (\ref{ThL2}), it suffices to show equality of the restrictions to $\Sl_{n-1}$ in the Grothendieck group. To prove this equality, by induction, we may assume 
the assertion to hold for $\Sl_{n-1}$, yielding 
\[
\begin{array}{rcl}
[S_{\sF_p}^\lambda]|_{\Sl_{n-1}}
& \aufgl{\bf\cite[\rm 9.2]{J}}  & \sum_{i\in R(\lambda)} [S_{\sF_p}^{\lambda(i)}] \\
& \aufgl{induction} & \sum_{i\in R(\lambda)} \Big( [D_{\sF_p}^{\lambda(i)}] + \{\mb{$i\neq 1$ and $\lambda_i\neq t$}\} [D_{\sF_p}^{\lambda(i)[t]}] \\
&                   &\hspace*{28.1mm} +\,\{\mb{$i = 1$ and $\lambda_{s+1} + 1 < t$}\} [D_{\sF_p}^{\lambda(1)[t-1]}] \Big)\; . \\
\end{array}
\]
On the other hand, since all removable nodes, except for $(s,t)$ if $\lambda_s = t$, are normal, and, using $p\, |\, h_t$, all normal nodes are good, 
modular branching {\bf\cite[\rm 0.6]{Kl}} enables us to calculate
\[
[D_{\sF_p}^\lambda]|_{\Sl_{n-1}}\; = \;\sum_{i\in R(\lambda)} \{\lambda_i\neq t\} [D_{\sF_p}^{\lambda(i)}]\; ,
\]
as well as
\[
[D_{\sF_p}^{\lambda[t]}]|_{\Sl_{n-1}}\;= \;\{\lambda_{s+1} + 1 < t\}[D_{\sF_p}^{\lambda(s)[t-1]}] 
+ \sum_{i\in R(\lambda)\ohne\{ s\}} [D_{\sF_p}^{\lambda(i)[t]}]\; .
\]
We are done by remarking that if $\lambda_s > t$, then $\lambda(s)[t] = \lambda(1)[t]$; that if $\lambda_s = t$, then $\lambda(s) = \lambda(1)[t]$;
and that if $\lambda_{s+1} + 1 < t$, then $\lambda(s)[t-1] = \lambda(1)[t-1]$. 

{\it Consider the case $\lambda_1 = \lambda_2$.} 
From $p\, |\, h_t$ we conclude that $\lambda_2 + \lambda'_1 \geq p+1$, and from $p > n - \lambda_1$ we infer that 
$\lambda_2 + \lambda'_1 < p+2$. Now, $\lambda_2 + \lambda'_1 = p + 1$ leads to $\lambda = (k^2,1^l)$, with $k\geq 1$, $l\geq 0$ and $k + l = p - 1$, so in 
particular $t = 1$.

If $k = 1$, then $[S_{\sF_p}^\lambda] = [D_{\sF_p}^\lambda] + [D_{\sF_p}^{\lambda[t]}]$, since $[D_{\sF_p}^{(1^p)}] = 0$ by convention, and since 
$[S_{\sF_p}^{(1^p)}] = [D_{\sF_p}^{(2,1^{p-2})}]$ by {\bf\cite[\rm Th.\ A]{JII}}.

If $k > 1$, we obtain one the one hand, by induction,
\[
[S_{\sF_p}^\lambda]|_{\Sl_{n-1}}\; =\; [D_{\sF_p}^{(k,k-1,1^l)}] + \{ l\geq 1\} \left([D_{\sF_p}^{(k+1,k-1,1^{l-1})}] + [D_{\sF_p}^{(k^2,1^{l-1})}]\right) 
+ \{ l = 0\}[D_{\sF_p}^{(n-1)}] \; .
\]
On the other hand, modular branching {\bf\cite[\rm 0.6]{Kl}} yields
\[
\begin{array}{rcl}
[D_{\sF_p}^{(k^2,1^l)}]|_{\Sl_{n-1}}                  & = & [D_{\sF_p}^{(k,k-1,1^l)}] \\
\{l\geq 1\}[D_{\sF_p}^{(k+1,k,1^{l-1})}]|_{\Sl_{n-1}} & = & \{l\geq 1\}\left([D_{\sF_p}^{(k+1,k-1,1^{l-1})}] + [D_{\sF_p}^{(k^2,1^{l-1})}]\right) \\
\{ l = 0\}{[D_{\sF_p}^{(n)}]}|_{\Sl_{n-1}}            & = & \{ l = 0\}[D_{\sF_p}^{(n-1)}]\; . \\
\end{array}
\]
Thus
\[
\begin{array}{rcl}
[S_{\sF_p}^\lambda] 
& = & [S_{\sF_p}^{(k^2,1^l)}] \\
& = & [D_{\sF_p}^{(k^2,1^l)}] + \{l\geq 1\}[D_{\sF_p}^{(k+1,k,1^{l-1})}] + \{ l = 0\}[D_{\sF_p}^{(n)}] \\ 
& = & [D_{\sF_p}^\lambda] + [D_{\sF_p}^{\lambda[t]}]\; .\\ 
\end{array}
\]
\end{Proposition}

\begin{Theorem}
\label{ThL4}
Let $\lambda$ be a partition of $n$, let $p$ be a prime such that $p > n - \lambda_1$. Let $h_j$ and, if existent, $t$ be as in (\ref{NotL1_5}).

If $p$ does not divide $h_j$ for any $j\in [1,\lambda_2]$, then
\[
[S_{\sF_p}^\lambda] \; =\; [D_{\sF_p}^\lambda]_0\; .
\]
If $p$ divides $h_t$, then
\[
[S_{\sF_p}^\lambda] \; =\; [D_{\sF_p}^\lambda]_0 + [D_{\sF_p}^{\lambda[t]}]_{v_p(h_t)}\; .
\]

\rm
If $p$ divides $h_t$, Schaper's formula (see e.g.\ {\bf\cite[\rm 5.33]{Mathas}}) reads
\[
[S_{\sZ_{(p)}}^{\lambda,\ast}/S_{\sZ_{(p)}}^\lambda]\; =\; \sumd{i\in [2,s]} (-1)^{s-i} v_p(h_t) [S_{\sF_p}^{\lambda\spi{i}}]\; ,
\]
where the diagram of the partition of $\lambda\spi{i}$ is obtained by removing the skew hook belonging to the node $(i,t)$ from the diagram of $\lambda$, and
by attaching as many nodes to the first row as there are in this skew hook (cf.\ {\bf\cite[\rm p.\ 223]{JM}}). 

Now, if $h\spi{i}_j$ denotes the hook length of the node $(1,j)$ of the diagram of $\lambda\spi{i}$, where $j\in [1,\lambda_2]$, and if $p$ divides 
$h\spi{i}_{t\spi{i}}$, then $\lambda\spi{i}[t\spi{i}] = \lambda\spi{i-1}$ for $i\in [3,s]$. Moreover, $\lambda\spi{s} = \lambda[t]$, and 
$S_{\sF_p}^{\lambda\spi{2}} = D_{\sF_p}^{\lambda\spi{2}}$ by (\ref{ThL1}). Thus by (\ref{PropL3}), we have
\[
[S_{\sZ_{(p)}}^{\lambda,\ast}/S_{\sZ_{(p)}}^\lambda]\; =\; v_p(h_t) [D_{\sF_p}^{\lambda[t]}]\; ,
\]
and we are done by (\ref{CorN2}). 
\end{Theorem}

\begin{Corollary}
\label{CorL6}
Suppose $p$ divides $h_t$. As modules over $\Z_{(p)}$, we have
\[
S_{\sZ_{(p)}}^{\lambda,\ast}/S_{\sZ_{(p)}}^\lambda\;\iso\; \left(\Z_{(p)}/h_t\Z_{(p)}\right)^{\dim_{\ssF_p} D_{\ssF_p}^{\lambda[t]}}\; ,
\]
where, using the notation of the proof of (\ref{ThL4}), we have
\[
\dim_{\sF_p} D_{\sF_p}^{\lambda[t]}\; =\; \sum_{i\in [2,s]} (-1)^{s-i} \dim_{\sF_p} S_{\sF_p}^{\lambda\spi{i}}\; .
\]

\rm
This follows using (\ref{ThL4}), its proof and (\ref{RemN0}).
\end{Corollary}

\section{Transposition}

Let $n\geq 1$. Given a $\Z\Sl_n$-module $X$, we denote $X^- := X\ts_\sZ S^{(1^n)}$. Given a partition $\lambda$, we denote $n_\lambda := \rk_\sZ\, S^\lambda$. 

\begin{Proposition}
\label{PropTr1}
Suppose given $n\geq 1$, $p$ a prime and $\lambda$ a partition of $n$. For $i\in [0,v_p(n!/n_\lambda)]$, we have
\[
\frac{S_{\sF_p}^\lambda(i)}{S_{\sF_p}^\lambda(i + 1)}
\;\iso\; \left(\frac{S_{\sF_p}^{\lambda'}(v_p(n!/n_\lambda) - i)}{S_{\sF_p}^{\lambda'}(v_p(n!/n_\lambda) - i + 1)}\right)^{\!\! -}\; .
\]
For $i > v_p(n!/n_\lambda)$, we have $S_{\sF_p}^\lambda(i)/S_{\sF_p}^\lambda(i + 1) = 0$.

\rm
By {\bf\cite[\rm 6.7]{J}} (and {\bf\cite[\rm 6.2.5]{K99}}), we have an isomorphism
\[
\begin{array}{rcl}
S^\lambda & \lraiso  & S^{\lambda',\ast,-} \\
\spi{a}   & \lramaps & (\{ a'\}, -)\; . \\
\end{array}
\]
Given a $\lambda$-tableau $[a]$, we let $C_a$ denote the column stabilizer of $[a]$ in $\Sl_n$, and $\kappa_a^- := \sum_{\sigma\in C_a} \sigma\eps_\sigma$; we let
$R_a$ denote the row stabilizer of $[a]$ in $\Sl_n$ and $\rho_a^+ := \sum_{\sigma\in R_a} \sigma$. Consider the 
composition
\[
\begin{array}{rcccccccl}
S^\lambda & \lraa{\eta}    & S^{\lambda,\ast}     & \lraiso  & S^{\lambda',-}     & \lraa{\eta^-} & S^{\lambda',\ast,-}           & \lraiso  & S^\lambda \\
\spi{a}   &\lramaps        & (\{ a\},-)\kappa^-_a & \lramaps & \spi{a'}\kappa^-_a & \lramaps      & (\{ a'\},-)\rho^+_a\kappa^-_a & \lramaps & 
                                                                                               \spi{a}\rho_a^+\kappa^-_a \; =\;\fracd{n!}{n_\lambda}\spi{a}\; ,\\
\end{array}
\]
the last equality following from {\bf\cite[\rm 23.2 ii]{J}}. We localize at $p$ and apply (\ref{RemR1}) to $X = S_{\sZ_{(p)}}^\lambda$, 
$Y = S_{\sZ_{(p)}}^{\lambda',-}$ and $N = v_p(n!/n_\lambda)$.
\end{Proposition}

\begin{Corollary}
\label{CorTr2}
We have $S_{\sF_p}^\lambda(v_p(n!/n_\lambda))/S_{\sF_p}^\lambda(v_p(n!/n_\lambda) + 1) \iso D_{\sF_p}^{\lambda',-}$.
\end{Corollary}

\begin{Corollary}
\label{CorTr3}
For $i\in [1,n_\lambda]$, the product of the $i$th elementary divisor of the Gram matrix of $S^\lambda$ and the $(n_\lambda + 1 - i)$th elementary 
divisor of the Gram matrix of $S^{\lambda'}$ yields $n!/n_\lambda$. 

In particular, the largest elementary divisor of the Gram matrix of $S^\lambda$ divides $n!/n_\lambda$. 

\rm
The first argument is to localize at a prime $p$ and to apply (\ref{RemN0}) to (\ref{PropTr1}). In this sense, (\ref{PropTr1}) is the module version of 
(\ref{CorTr3}).

The second argument is to consider the composition used in the proof of (\ref{PropTr1}) directly, which yields the result without reverting to localization.
\end{Corollary}

\begin{footnotesize}
\begin{quote}
\begin{Remark}\rm
A Brauer-Nesbitt type argument (cf.\ section \ref{SubSecQuasi}) shows that if $G$ is a finite group, $R$ a discrete valuation ring of characteristic zero with 
fraction field $K$ splitting $G$, and $X$ a simple $RG$-lattice carrying a nondegenerate $G$-invariant $R$-bilinear form, then the quotient of the last 
elementary divisor by the first elementary divisor of the Gram matrix of this bilinear form divides $|G|/\rk_R X$ in $R$. In this generality, there is no 
substitute for transposition of partitions. Nonetheless, in our particular case this fact can be derived using transposition of partitions (\ref{CorTr3}).
\end{Remark}
\end{quote}
\end{footnotesize}

The following has been discovered by {\sc C.\ D.\ Gay} and rederived by {\sc James} and {\sc Murphy} (cf.\ {\bf\cite[\rm p.\ 234]{JM}}).

\begin{Corollary}
\label{CorTr4}
The product of the determinants of the Gram matrices of $S^\lambda$ and of $S^{\lambda'}$ is given by $(n!/n_\lambda)^{n_\lambda}$.
\end{Corollary}

For later use, we record the related

\begin{Lemma}
\label{LemTr2}
We have $(\spi{a}\rho_a^+,-) = \fracd{n!}{n_\lambda}(\{a \},-)$. 

\rm 
To this end, we compose the string of morphisms displayed in the proof of (\ref{PropTr1}) from $S^{\lambda,\ast}$ to $S^{\lambda,\ast}$. 
\end{Lemma}

Sometimes, we shall abbreviate $\spi{a}\rho_a^+ =: \spi{a}\rho^+$.

\section{Explicit bases}

\begin{footnotesize}
\begin{quote}
So far, we have been calculating elementary divisors of Gram matrices for Specht modules without reference to any explicit diagonalization. Instead, we made use 
of known decomposition numbers as well as of Schaper's formula.

It would be preferable to construct two bases diagonalizing the Gram matrix, but this seems to be difficult in general. For hook partitions, however, explicit 
diagonalization is a simple (and precise) way to get the elementary divisors (\ref{ThHook}). We also treat the partition $(2^2,1^{n-4})$
in this manner (\ref{Th2^2}), but already this modest case indicates how complicated such bases seem to look like in general, at least in terms of polytabloids. 
We remark that even for simple modules whose dimension is known (such as two-row partitions, three-row partitions with third row $< p$, or those treated in 
{\bf\cite{Mathieu}}), explicit bases are largely unknown (or at least unpublished). Such a basis, lifted to the Specht module of the same partition, could 
locally be used as part of a diagonalizing basis. The bases found so far diagonalize globally.

To simplify, we spare ourselves the construction of the second basis by giving the first basis in a convenient manner with respect to the standard choice of the 
second, see section \ref{SubSecLA}. By a matrix inversion of a unipotent upper triangular integral matrix, a diagonalizing second basis ensues.
\end{quote}
\end{footnotesize}

\subsection{A linear algebra lemma}
\label{SubSecLA}
\begin{Lemma}
\label{LemFund}
Let $X$ be a finitely generated free $\Z$-module of rank $m\geq 1$, let $(-,=)$ be a nondegenerate bilinear form on $X$. Let $G$ be the Gram matrix of that 
bilinear form with respect to some $\Z$-linear basis. Let $(x_1,\dots,x_m)$, $(y_1,\dots,y_m)$ be two tuples of elements of $X$ such that the following hold.
\begin{itemize}
\item[(i)] We have $(x_i,y_i)\neq 0$ for all $i\in [1,m]$.
\item[(ii)] We have $(x_i,y_j) = 0$ for all $i,j\in [1,m]$ such that $i > j$.
\item[(iii)] $(x_i,y_i)$ divides $(x_i,y_j)$ for all $i,j\in [1,m]$ such that $i < j$. 
\item[(iv)] $(x_i,y_i)$ divides $(x_{i+1},y_{i+1})$ for all $i\in [1,m-1]$.
\item[(v)] $\prod_{i\in [1,m]} (x_i,y_i) = \pm\det G$.
\end{itemize}
Then the tuples $(x_1,\dots,x_m)$ and $(y_1,\dots,y_m)$ are $\Z$-linear bases of $X$. The elementary divisors of $G$ are given by the tuple of integers
\[
\Big((x_1,y_1), (x_2,y_2), \dots, (x_n,y_n)\Big)\; .
\]

\rm
Let $g:X\lra X^\ast$ send $\xi\in X$ to $(\xi,-)$. 
The tuple $(x_1,\dots,x_m)$ is linear independent over $\Z$, as its image in $X^\ast$ via $g$ shows, using (i,\,ii). Likewise, the tuple $(y_1,\dots,y_m)$ is 
linearly independent. In the composite embedding
\[
\Z\spi{x_1, \dots, x_m}\;\hraa{f}\; X\;\hraa{g}\; X^\ast\;\hraa{h}\; \Z\spi{y_1, \dots, y_m}^\ast\;,
\]
both $f\!gh$ and $g$ have determinant $\pm\det G$, by (v). Thus $f$ and $h$ are equalities. Now (i,\,ii,\,iii,\,iv) ensure that the elementary divisors result 
as claimed.
\end{Lemma}

The element $y_i$ of the second basis is called the {\it diagonal correspondent} of the element $x_i$ of the first basis.

\subsection{Bases for hooks}

\begin{footnotesize}
\begin{quote}
The result of this section has independently been observed by {\sc Mathas} and {\sc James} [unpublished]. As ingredient for (\ref{LemFund} v), we shall make use 
of the determinant of the Gram matrix, but not of Schaper's formula itself.
\end{quote}
\end{footnotesize}

Let $n\geq 2$, let $l\in [0,n-1]$, and let $\lambda := (n-l,1^l)$ be a {\it hook partition} of $n$. Note that a $\lambda$-polytabloid is 
determined by the tuple of entries in the first column, read from top to bottom. So given a tuple $b = (b_0,b_1,\dots, b_l)$, the entries $b_i\in [1,n]$ 
pairwise distinct, we shall write $\spi{b}$ for the $\lambda$-polytabloid it determines. In the sequel, by a tuple we understand a tuple with pairwise distinct 
entries. To the underlying set of a tuple we refer without further comment.

Counting ordered tuples, we obtain $\rk_\sZ\, S^\lambda = \smatze{n-1}{l}$.

If $b$ is a tuple and $x$ an element, then $(x,b)$ denotes the tuple $b$ with an entry $x$ appended on the left, $(b,x)$ denotes the tuple $b$ with an entry $x$ 
appended on the right, etc. If tuple brackets appear within tabloid or polytabloid brackets, we omit the tuple brackets. For instance, if $n = 5$, $l = 2$ and 
$b = (5,3)$, then $\spi{2,b} = \smatspdd{2}{1}{4}{5}{}{}{3}{}{}$.

\begin{Lemma}
\label{LemHook3}
We have
$\;
\# (S^{\lambda,\ast}/S^\lambda)\; =\; l!^{\smatze{n-1}{l}}\cdot n^{\smatze{n-2}{l-1}}\; . 
$

\rm
This follows by induction on $n$ using the Branching Theorem for Determinants {\bf\cite[\rm p.\ 225]{JM}}. 
\end{Lemma}

\begin{Lemma}
\label{LemHook1}
Suppose given increasingly ordered tuples $b,c\tm [2,n-1]$ of length $\# b = \# c = l$. We have
\[
(\spi{n,b},\spi{1,c})_0 \; =\; 
\left\{
\begin{array}{lll}
1          & \mb{if} & b = c\; , \\
0          & \mb{if} & b\neq c\; .  \\
\end{array}
\right.
\]
\end{Lemma}

\begin{Lemma}
\label{LemHook2}
Suppose given increasingly ordered tuples $d\tm [2,n-1]$ of length $\# d = l - 1$, and $c\tm [2,n]$ of length $\# c = l$. We have
\[
\Big(\sumd{s\in [1,n-1]\ohne d}\spi{s,d,n},\spi{1,c}\Big)_0 \; =\; 
\left\{
\begin{array}{lll}
n          & \mb{if} & (d,n) = c\; , \\
0          & \mb{if} & (d,n)\neq c\; . \\
\end{array}
\right.
\]

\rm
In fact,
\[
\begin{array}{rcl}
\Big(\sumd{s\in [1,n-1]\ohne d}\spi{s,d,n},-\Big)_0 
& = & (n-l-1)!^{-1}\Big(\spi{1,d,n}\rho^+, -\Big)_0 \\
&\aufgl{(\ref{LemTr2})} & n\cdot\Big(\{ 1,d,n\}, -\Big)\; . \\
\end{array}
\]
\end{Lemma}

\begin{Theorem}
\label{ThHook}
As abelian groups, we have
\[
S^{\lambda,\ast}/S^\lambda\; =\; S^{(n-l,1^l),\ast}/S^{(n-l,1^l)}\;\iso\; (\Z/l!\,\Z)^{\smatze{n-2}{l}}\ds (\Z/n\cdot l!\,\Z)^{\smatze{n-2}{l-1}}\; .
\]
Bases $x$ and $y$ that trigonalize the Gram matrix in the sense of (\ref{LemFund}) are given by a tuple $y$ consisting of the standard basis, ordered in such a 
manner that all elements without entry $n$ in the first column are before all elements with entry $n$ in the first column, and by a tuple $x = (x',x'')$, with 
$x'$ consisting of the elements of the form $\spi{n,b}$ as in (\ref{LemHook1}), with respective diagonal correspondent $\spi{1,b}$ in $y$, and $x''$ consisting 
of elements of the form $\sum_{s\in [1,n]\ohne d}\spi{s,d,n}$ as in (\ref{LemHook2}), with respective diagonal correspondent $\spi{1,d,n}$ in $y$. 
\end{Theorem}

\begin{footnotesize}
\begin{quote}
\begin{Example}
\label{ExHook4}\rm
Let $n = 4$ and $\lambda = (2,1^2)$. With respect to 
\[
x \; =\; \left(\smatspdz{4}{1}{2}{}{3}{},\smatspdz{1}{3}{2}{}{4}{} + \smatspdz{3}{1}{2}{}{4}{},
               \smatspdz{1}{2}{3}{}{4}{} + \smatspdz{2}{1}{3}{}{4}{}\right)\; ,\;\;\; 
y \; =\; \left(\smatspdz{1}{4}{2}{}{3}{},\smatspdz{1}{3}{2}{}{4}{},\smatspdz{1}{2}{3}{}{4}{}\right)\; , 
\]
the Gram matrix $\left(\;(x_i,y_j)\;\right)_{i,j}$ takes the form $\smatdd{2}{-2}{2}{0}{8}{0}{0}{0}{8}$.
\end{Example}
\end{quote}
\end{footnotesize}

\begin{Remark}[Jantzen subquotients]
\label{RemHook5}\rm
Let $p$ be a prime. To consider the Jantzen subquotients, we have to distinguish two cases. 

{\bf Case $p\geq 3$.} Given a partition $\nu$, we denote by $\nu^\mb{\scr r}$ the $p$-regularized partition in the sense of {\bf\cite[\rm p.\ 46]{JII}}, where 
the corresponding diagram construction is explained -- the {\it diagram} of $\nu$ being given by 
$[\nu] := \{ (i,j)\in\Z_{\geq 1}\ti\Z_{\geq 1}\; |\; j\leq\nu_i\}$. 

Now, if $p$ does not divide $n$, then $S_{\sF_p}^\lambda \iso D_{\sF_p}^{\lambda^\mb{\sscr r}}$ is simple {\bf\cite[\rm p.\ 52]{JII}}. By (\ref{ThHook}), we 
obtain
\[
[S_{\sF_p}^\lambda] \; =\; [D_{\sF_p}^{\lambda^\mb{\sscr r}}]_{v_p(l!)}\; .
\]
So suppose that $p$ divides $n$. 

If $l = 0$, we obtain 
\[
[S_{\sF_p}^\lambda] \; =\;  [S_{\sF_p}^{(n)}] \; =\;  [D_{\sF_p}^{(n)}]_0\; . 
\]
If $l = n-1$, we obtain by {\bf\cite[\rm Th.\ A]{JII}} 
\[
[S_{\sF_p}^\lambda] \; =\; [S_{\sF_p}^{(1^n)}] \; =\; [D_{\sF_p}^{(1^n)^\mb{\sscr r}}]_{v_p(n!)}\; . 
\]
If $l\in [1,n-2]$, the decomposition numbers of {\bf\cite[\rm p.\ 52]{JII}} and the long exact hook sequence of {\bf\cite[\rm Lem.\ 2]{Pe}} 
(cf.\ {\bf\cite[\rm 4.2.3, 4.2.4]{K99}}) give
\[
[S_{\sF_p}^\lambda]\; =\; [S_{\sF_p}^{(n-l,1^l)}]\; =\; [D_{\sF_p}^{(n - l - \{\, l\;\geq\; n - n/p\},1^\ast)^\mb{\sscr r}}]_{v_p(l!)} 
+ [D_{\sF_p}^{(n - (l-1) - \{\, (l-1)\;\geq\; n - n/p\},1^\ast)^\mb{\sscr r}}]_{v_p(l!\, n)}\; , 
\]
where we abbreviated $(n-j,1^\ast) := (n-j,1^j)$. In fact, the image in $S_{\sF_p}^\lambda$ of the differential of the long exact hook sequence equals 
$S_{\sF_p}^\lambda(l!\, n)$, as results from a comparison of {\bf\cite[\rm 4.2.3]{K99}} with (\ref{LemHook2}).

{\bf Case $p = 2$.} Suppose given $k,m\geq 0$, written $2$-adically as $k + 1 = \sum_{s\in [0,K]} a_s 2^s$ and $m = \sum_{s\in [0,M]} b_s 2^s$, where 
$a_s,b_s\in\{ 0,1\}$, where $a_K = 1$ and where $b_M = 1$ if $m\geq 1$, and $M = -1$ if $m = 0$. Assume $K > M$. Let 
$t := \max\left(\{ s\in [1,M]\; |\; a_s < b_s\}\cup \{ 0\}\right) $, and define 
\[
F_2(k,m) \; :=\; \sum_{i\geq 0} f_2(k,m-2i) = 2^{\sum_{u\in [1,t]} a_u}
                                              \cdot\left(1 + \!\!\sum_{s\in [t+1,M]}\!\! a_s b_s\cdot 2^{\sum_{u\in [t+1,s-1]} a_u}\right)\cdot\{ km\con_2 0\}\; ,
\]
the notation being as for (\ref{ThT2}).

First, suppose $2l\leq n-1$. For $j\in [0,l]$, {\bf\cite[\rm p.\ 93]{J}} and {\bf\cite[\rm 24.15]{J}} yield 
$[S_{\sF_2}^\lambda : D_{\sF_2}^{(n-j,j)}] = F_2(n-2j,l-j)$ (cf.\ {\bf\cite{Su}}). If $n$ is even, we use the same argument as in the case $p\geq 3$, if $n$ is 
odd, there is no further argument necessary, to conclude 
\[
[S_{\sF_2}^\lambda] \; = \; [S_{\sF_2}^{(n-l,1^l)}] \; =\; \sum_{j\in [0,l]} F_2(n-2j,l-j)\, [D_{\sF_2}^{(n-j,j)}]_{v_2(l!\, n^{\{ l \not\con_2 j\}})}\; .
\]
Next, suppose $2l \geq n-1$. Since $[S_{\sF_2}^\lambda] = [S_{\sF_2}^{\lambda',\ast,-}] = [S_{\sF_2}^{\lambda'}]$, we obtain in the same manner
\[
[S_{\sF_2}^\lambda] \; = \; [S_{\sF_2}^{(n-l,1^l)}]\; = \; [S_{\sF_2}^{(l+1,1^{n-l-1})}] 
\; = \sum_{j\in [0,n-l-1]} F_2(n-2j,n-l-1-j)\, [D_{\sF_2}^{(n-j,j)}]_{v_2(l!\, n^{\{ l \not\con_2 j\}})}\; .
\]
\end{Remark}

\begin{quote}
\begin{footnotesize}
\begin{Example}
\label{ExHook6}\rm
We have 
\[
[S_{\sF_2}^{(6,1^8)}] \; =\;  3 [D_{\sF_2}^{(14,0)}]_7 + 2 [D_{\sF_2}^{(13,1)}]_8 + 2 [D_{\sF_2}^{(12,2)}]_7 + [D_{\sF_2}^{(11,3)}]_8 
+ [D_{\sF_2}^{(10,4)}]_7 + [D_{\sF_2}^{(9,5)}]_8\; .
\]
\end{Example}
\end{footnotesize}
\end{quote}

\subsection{Gram matrix entries for two-column partitions}
\label{SubSecTwoCol}

Let $n\geq 1$, $1\leq h\leq n/2$ and $\lambda = (2^h,1^{n-2h})$. Recall that tuples have pairwise distinct entries by convention.

\begin{Lemma}
\label{LemF}
Let $i\in [0,h]$. Suppose given pairwise disjoint tuples $\phi,\psi,\theta,\zeta\tm [1,n]$ of length
\[
\begin{array}{lcl}
\#\phi    & = & n - 2h + i\; , \\
\#\psi   & = & i\; , \\
\#\theta & = & h - i\; ,\\  
\#\zeta  & = & h - i\; . \\
\end{array}
\]
Let $[a]$ be the tableau with first column $(\theta,\phi)$ and second column $(\zeta,\psi)$. Let $[b]$ be the tableau with first column $(\zeta,\phi)$ and 
second column $(\theta,\psi)$. Then
\[
\Big(\spi{a}, \spi{b}\Big)\; =\; h!\cdot (h-i)!\cdot (n - 2h + i)!\; .
\]

\rm
We have to calculate
\unitlength0.07mm        
\begin{center}
\begin{picture}(1400,600)
\put( -50, 375){$\left(\rule[12mm]{0mm}{0mm}\right.$}
\put(   0, 375){\line(1,4){60}}
\put(   0, 375){\line(1,-4){60}}
\put( 100, 600){\line(1,0){70}}
\put( 100, 400){\line(1,0){70}}
\put( 100, 600){\line(0,-1){200}}
\put( 170, 600){\line(0,-1){200}}
\put( 122, 490){$\scm\theta$}
\put( 100, 350){\line(1,0){70}}
\put( 100,   0){\line(1,0){70}}
\put( 100, 350){\line(0,-1){350}}
\put( 170, 350){\line(0,-1){350}}
\put( 122, 165){$\scm\phi$}
\put( 200, 600){\line(1,0){70}}
\put( 200, 400){\line(1,0){70}}
\put( 200, 600){\line(0,-1){200}}
\put( 270, 600){\line(0,-1){200}}
\put( 222, 490){$\scm\zeta$}
\put( 200, 350){\line(1,0){70}}
\put( 200, 150){\line(1,0){70}}
\put( 200, 350){\line(0,-1){200}}
\put( 270, 350){\line(0,-1){200}}
\put( 222, 240){$\scm\psi$}
\put( 370, 375){\line(-1,4){60}}
\put( 370, 375){\line(-1,-4){60}}
\put( 390, 375){\large $,$} 
\put( 430, 375){\line(1,4){60}}
\put( 430, 375){\line(1,-4){60}}
\put( 530, 600){\line(1,0){70}}
\put( 530, 400){\line(1,0){70}}
\put( 530, 600){\line(0,-1){200}}
\put( 600, 600){\line(0,-1){200}}
\put( 552, 490){$\scm\zeta$}
\put( 530, 350){\line(1,0){70}}
\put( 530,   0){\line(1,0){70}}
\put( 530, 350){\line(0,-1){350}}
\put( 600, 350){\line(0,-1){350}}
\put( 552, 165){$\scm\phi$}
\put( 630, 600){\line(1,0){70}}
\put( 630, 400){\line(1,0){70}}
\put( 630, 600){\line(0,-1){200}}
\put( 700, 600){\line(0,-1){200}}
\put( 652, 490){$\scm\theta$}
\put( 630, 350){\line(1,0){70}}
\put( 630, 150){\line(1,0){70}}
\put( 630, 350){\line(0,-1){200}}
\put( 700, 350){\line(0,-1){200}}
\put( 652, 240){$\scm\psi$}
\put( 800, 375){\line(-1,4){60}}
\put( 800, 375){\line(-1,-4){60}}
\put( 825, 375){$\left.\rule[12mm]{0mm}{0mm}\right)$}
\put( 880, 375){$\; =\;\sumd{\tau\,\in\, C_a}\;\;\sumd{\sigma\,\in\, C_b}\; \Big\{\{ a\}\tau = \{ b\}\sigma\Big\}\eps_\tau\eps_\sigma\;$.}
\end{picture}
\end{center}
\unitlength0.1mm        
First, assume that $\zeta$ and $\psi$ remain fixed in $[a]$, i.e.\ that $\tau$ restricts to the identity on $\zeta\cup\psi$. In order to obtain 
$\{ a\}\tau = \{ b\}\sigma$, the permutation $\sigma$ has to restrict to the identity on $\zeta\cup\psi$ as well. Moreover, $\sigma$ and $\tau$ have to restrict 
to the same automorphism of $\theta$, whence the factor $(h-i)!$, and to the same automorphism on $\phi$, whence the factor $(n - 2h + i)!$. Next, removing the 
assumption, we remark that we may simultaneously permute the rows $1$ to $h$ of $\{ a\}\tau$ and of $\{ b\}\sigma$, whence the factor $h!$.
\end{Lemma}

\begin{Lemma}
\label{LemF2}
Suppose given increasingly ordered tuples $\xi,\eta\tm [1,n]$ of length $h$ with intersection of cardinality $i$. We write $\spi{\xi}$ for the polytabloid with 
second column $\xi$, and first column $[1,n]\ohne\xi$, increasingly ordered from top to bottom. Let $\sum\xi$ denote the sum of the entries of $\xi$. Then
\[
(\spi{\xi},\spi{\eta})_0\; =\; \fracd{(h-i)!(n-2h + i)!}{(n-2h)!}\cdot (-1)^{h - i + \sum\xi + \sum\eta}\ .
\]

\rm
This follows from (\ref{LemF}) by reordering the columns of the polytabloids involved.
\end{Lemma}

\subsection{Bases for $(2^2,1^{n-4})$}
\label{SubSec221n-4}

Let $n\geq 6$, let $\lambda = (2^2,1^{n-4})$. Given $a,b\in [1,n]$, $a\neq b$, we denote by $\smatspze{a}{b}$ the $\lambda$-polytabloid with second column 
$(a,b)$, and whose first column is increasingly ordered from top to bottom. That is, we use the notation that has already been employed in (\ref{LemF2}).

\begin{Lemma}
\label{Lem2^2_1}
We have
\[
\#(S^{\lambda,\ast}/S^\lambda) \; =\; (2(n-4)!)^{(n^2-3n)/2}\cdot (n-1)^{(n^2 - 3n - 2)/2}\cdot (n-2)^{(n^2 - 5n + 2)/2}\cdot 2\; .
\]

\rm
The transposed partition has $\#(S^{(n-2,2),\ast}/S^{(n-2,2)}) = (n-1)\cdot (n-2)^{n-1}/2$, as we take from the first example on {\bf\cite[\rm p.\ 224]{JM}}. 
Since $\rk_\sZ\, S^\lambda = n(n-3)/2$ by {\bf\cite[\rm 20.1]{J}}, the assertion ensues from (\ref{CorTr4}).
\end{Lemma}

\begin{Lemma}
\label{Lem2^2_2}
We obtain the following table for values of $(-,=)_0 = \frac{1}{2(n-4)!}(-,=)$, in which the columns are indexed by the standard polytabloid basis of $S^\lambda$.

{\rm
\begin{footnotesize}
\begin{tabular}{l||l|l|l|}
$(-,=)_0$ & $\smatspze{b}{n}$ & $\smatspze{b}{n-1}$ & $\smatspze{b}{c}$          \\
          & $b\in [2,n-1]$    & $b\in [2,n-2]$      & $b\in [2,n-2]$,          \\
          &                   &                     & $c\in [4,n-2]$, $b < c$ \\\hline\hline
$\smatspze{3}{n}$
&
$\begin{array}{l}
(n-2)(n-3)\\ 
\mb{ if $b = 3$}       \\
(-1)^b\cdot (n-3)     \\
\mb{ if $b\neq 3$}     \\
\end{array}$
&
$\begin{array}{l}
(n-3)                  \\
\mb{ if $b = 3$}       \\
(-1)^b\cdot 2          \\
\mb{ if $b\neq 3$}     \\
\end{array}$
&
$\begin{array}{l}
(-1)^{c+n+1}\cdot (n-3) \\ 
\mb{ if $b = 3$}         \\
(-1)^{b+c+n+1}\cdot 2   \\
\mb{ if $b\neq 3$}       \\
\end{array}$
\\\hline
$\smatspze{3}{4}$
&
$\begin{array}{l}
(-1)^{n+1}\cdot (n-3)        \\ 
\mb{ if $b = 3$}              \\
(-1)^n\cdot (n-3)            \\ 
\mb{ if $b = 4$}              \\
(-1)^{b+n+1}\cdot 2          \\
\mb{ if $b\not\in\{ 3, 4\}$}  \\
\end{array}$
&
$\begin{array}{l}
(-1)^n\cdot (n-3)        \\ 
\mb{ if $b = 3$}              \\
(-1)^{n+1}\cdot (n-3)            \\ 
\mb{ if $b = 4$}              \\
(-1)^{b+n}\cdot 2                 \\
\mb{ if $b\not\in\{ 3, 4\}$}  \\
\end{array}$
&
$\begin{array}{l}
(n-2)(n-3)                          \\ 
\mb{ if $b = 3$, $c = 4$}            \\
(n-3)                               \\ 
\mb{ if $b = 2$, $c = 4$}            \\
(-1)^{c+1}\cdot (n-3)               \\ 
\mb{ if $b = 3$, $c\neq 4$} \\
(-1)^c\cdot  (n-3)                  \\ 
\mb{ if $b = 4$, $c\neq 4$} \\
(-1)^{b+c+1}\cdot 2                 \\ 
\mb{ if $b,c\not\in\{ 3,4\}$}        \\
\end{array}$
\\\hline
$
\begin{array}{l}
\smatspze{k}{n-1} + (-1)^k\smatspze{1}{n-1} \\
k\in [2,n-2]                                \\
\end{array}
$
&
$\begin{array}{l}
(n-1)                        \\
\mb{ if $b = k$}              \\
0                            \\
\mb{ if $b \neq k$}           \\
\end{array}$
&
$\begin{array}{l}
(n-1)(n-3)                   \\
\mb{ if $b = k$}              \\
0                            \\
\mb{ if $b \neq k$}           \\
\end{array}$
&
$\begin{array}{l}
(-1)^{n+c}\cdot (n-1)               \\
\mb{if $k = b$}                     \\
(-1)^{n+b}\cdot (n-1)               \\
\mb{if $k = c$}                     \\
0                                   \\
\mb{if $k\not\in\{ b,c\}$}          \\
\end{array}$
\\\hline
$
\begin{array}{l}
\smatspze{3}{n}            \\
+ (-1)^n \smatspze{n-1}{n} \\
- (n-3)\smatspze{3}{n-1}   \\
\end{array}
$
& 
$0$            
& 
$\begin{array}{l}
-(n-1)(n-3)(n-4)           \\
\mb{ if $b = 3$}            \\
(-1)^{b+1}\cdot (n-1)(n-4) \\ 
\mb{ if $b\neq 3$}          \\
\end{array}$
&
$\begin{array}{l}
(-1)^{c+n+1}\cdot (n-1)(n-4) \\ 
\mb{ if $b = 3$}              \\
(-1)^{b+c+n+1}\cdot 2(n-1)   \\
\mb{ if $b\neq 3$}            \\
\end{array}$
\\\hline
$
\begin{array}{l}
\smatspze{k}{l}\rho^+               \\
k\in [2,n-1],                    \\ 
l\in [4,n-1], k < l \\
\end{array}
$
&
0
& 
$\begin{array}{l}
(n-1)(n-2)                          \\
\mb{ if $k = b$, $l = n-1$}         \\
0                                   \\
\mb{ else}                          \\
\end{array}$
&
$\begin{array}{l}
(n-1)(n-2)                          \\
\mb{ if $k = b$, $l = c$}            \\
(-1)^{l+1}(n-1)(n-2)                \\
\mb{ if $k = c$, $b = 2$}            \\
0                                   \\
\mb{ else}                          \\
\end{array}$
\\\hline
\end{tabular}
\end{footnotesize}
}

\rm 
This table is to be verified using (\ref{LemF2}), except for the last row, for which we may most conveniently use (\ref{LemTr2}). Explicitely, we have for 
example $\smatspze{3}{4}\rho^+ = \smatspze{3}{4} - \smatspze{1}{4} + \smatspze{3}{2} + \smatspze{1}{2}$. 
\end{Lemma}

\begin{quote}
\begin{footnotesize}
\begin{Example}
\label{Ex2^2_3}
\rm
If $n = 6$, we obtain the following table. 

\begin{center}
\begin{tabular}{l||r|r|r|r|r|r|r|r|r|}
$(-,=)_0$         & $\!\!\smatspze{5}{6}\!\!$ & $\!\!\smatspze{4}{6}\!\!$ & $\!\!\smatspze{3}{6}\!\!$ & $\!\!\smatspze{2}{6}\!\!$ & $\!\!\smatspze{4}{5}\!\!$ 
& $\!\!\smatspze{3}{5}\!\!$ & $\!\!\smatspze{2}{5}\!\!$ & $\!\!\smatspze{3}{4}\!\!$ & $\!\!\smatspze{2}{4}\!\!$ \\\hline\hline
$\smatspze{3}{6}$ & $-3$ & $3$ & $12$ & $3$  & $2$  & $3$ & $2$ & $-3$ & $-2$ \\\hline
$\smatspze{3}{4}$ & $2$  & $3$ & $-3$ & $-2$ & $-3$ & $3$ & $2$ & $12$ & $3$  \\\hline\hline
$\smatspze{4}{5} + \smatspze{1}{5}$ & $0$ & $5$ & $0$ & $0$ & $15$ & $0$  & $0$  & $-5$ & $5$ \\\hline
$\smatspze{3}{5} - \smatspze{1}{5}$ & $0$ & $0$ & $5$ & $0$ & $0$  & $15$ & $0$  & $5$  & $0$ \\\hline
$\smatspze{2}{5} + \smatspze{1}{5}$ & $0$ & $0$ & $0$ & $5$ & $0$  & $0$  & $15$ & $0$  & $5$ \\\hline\hline
$\smatspze{3}{6} + \smatspze{5}{6} - 3\cdot\smatspze{3}{5}$ & $0$ & $0$ & $0$ & $0$ & $-10$ & $-30$ & $-10$ & $-10$ & $-10$ \\\hline
$\smatspze{4}{5}\rho^+$ & $0$ & $0$ & $0$ & $0$ & $20$ & $0$  & $0$  & $0$  & $20$ \\\hline\hline
$\smatspze{3}{5}\rho^+$ & $0$ & $0$ & $0$ & $0$ & $0$  & $20$ & $0$  & $0$  & $0$  \\\hline
$\smatspze{2}{5}\rho^+$ & $0$ & $0$ & $0$ & $0$ & $0$  & $0$  & $20$ & $0$  & $0$  \\\hline
$\smatspze{3}{4}\rho^+$ & $0$ & $0$ & $0$ & $0$ & $0$  & $0$  & $0$  & $20$ & $0$  \\\hline
$\smatspze{2}{4}\rho^+$ & $0$ & $0$ & $0$ & $0$ & $0$  & $0$  & $0$  & $0$  & $20$ \\\hline
\end{tabular}
\end{center}

Using (\ref{LemFund}, \ref{Lem2^2_1}), this yields 
\[
S^{(2^2,1^2),\ast}/S^{(2^2,1^2)}\;\iso\; (\Z/4\Z)^1\ds (\Z/20\Z)^3\ds (\Z/40\Z)^1\ds (\Z/80\Z)^4
\]
as abelian groups.
\end{Example}
\end{footnotesize}
\end{quote}

\begin{Theorem}[cf.\ (\ref{ExT7})]
\label{Th2^2}
Suppose $n\geq 6$. As abelian groups, we obtain
\[
\begin{array}{rcrl}
S^{(2^2,1^{n-4}),\ast}/S^{(2^2,1^{n-4})}
& \iso &     & (\Z/2^{\{ n\con_2 1\}}\cdot 2(n-4)!\,\Z)^1 \\
&      & \ds & (\Z/(n-1)\cdot 2(n-4)!\,\Z)^{n-3} \\
&      & \ds & (\Z/2^{\{n\con_2 0\}}(n-1)\cdot 2(n-4)!\,\Z)^1\\
&      & \ds & (\Z/(n-2)(n-1)\cdot 2(n-4)!\,\Z)^{(n^2 - 5n + 2)/2}\; . \\
\end{array}
\]
Trigonalizing bases $x$ and $y$ in the sense of (\ref{LemFund}) are given as follows. Let
\[
\begin{array}{rcll}
y_i            & := & \smatspze{n-i}{n}     & \mb{for $i\in [1,n-2]$} \\
y_{(n-2) + i}  & := & \smatspze{n-1-i}{n-1} & \mb{for $i\in [1,n-3]$} \\
y_{(2n-5) + i} & := & \smatspze{b_i}{c_i}   & \mb{for $i\in [1,\smatze{n-3}{2}\! - 1]$}\; , \\
\end{array}
\]
where $b_i,c_i\in [2,n-2]$, $b_i < c_i$, $(b_i, c_i)\neq (2,3)$, and where $i < j$ implies that $c_i > c_j$, or that $c_i = c_j$ and $b_i > b_j$ (i.e.\ we 
choose a reverse lexicographical ordering, read from bottom to top). Let
\[
\begin{array}{rcll}
x_i            & := & (-1)^{n-i}\smatspze{1}{n-1} + \smatspze{n-i}{n-1} & \mb{for $i\in [2,n-2]$} \\
x_{(n-2) + i}  & := & \smatspze{n-1-i}{n-1}\rho^+                       & \mb{for $i\in [2,n-3]$} \\
x_{(2n-5) + i} & := & \smatspze{b_i}{c_i}\rho^+                         & \mb{for $i\in [1,\smatze{n-3}{2} \! - 1]$}\; . \\
\end{array}
\]
Moreover, let
\[
x_1 := 
\left\{
\begin{array}{ll}
\smatspze{3}{4}                                & \mb{if $n$ is odd,} \\
\smatspze{3}{n} + \frac{n-2}{2}\smatspze{3}{4} & \mb{if $n$ is even,} \\
\end{array}
\right.
\]
and
\[
x_{(n-2) + 1} := 
\left\{
\begin{array}{rcrl}
\frac{n-3}{2}\smatspze{n-2}{n-1}\rho^+ & - & \frac{n-1}{2} \left(\smatspze{3}{n} + (-1)^n \smatspze{n-1}{n} - (n-3)\smatspze{3}{n-1}\right) & 
                                                                                                                                       \mb{if $n$ is odd,} \\
\smatspze{n-2}{n-1}\rho^+              & + & \left(\smatspze{3}{n} + (-1)^n \smatspze{n-1}{n} - (n-3)\smatspze{3}{n-1}\right)               & 
                                                                                                                                       \mb{if $n$ is even.} \\
\end{array}
\right.
\]

\rm
Conditions (i,\,ii,\,iii,\,iv) of (\ref{LemFund}) follow by (\ref{Lem2^2_2}), condition (v) follows by (\ref{Lem2^2_1}).
\end{Theorem}

\section{Miscellanea}

\subsection{Symmetric partitions}

A partition $\lambda$ is called {\it symmetric} if $\lambda = \lambda'$, otherwise, it is called {\it asymmetric.} We remark that by induction, ordinary 
branching {\bf\cite[\rm 9.2 ii]{J}} shows that if $\lambda\neq (1)$ is symmetric, then $n_\lambda := \rk_\sZ\, S^\lambda$ is even.

Suppose given a symmetric partition $\lambda\neq (1)$. Let $m$ be the {\rm middle jump factor,} i.e.\ the quotient of the $(n_\lambda/2 + 1)$st and the 
$(n_\lambda/2)$th elementary divisor of the Gram matrix of $S^\lambda$. Let $H := \prod_{i\in [1,s]} (2\lambda_i - 2i + 1)$ be the product of the main diagonal 
hook lengths, where $s = \lambda_s$.

\begin{Proposition}
\label{PropSym3}
The quotient $H/m$ is a square in $\Q$.

\rm
By (\ref{CorTr3}), the quotient $(n!/n_\lambda)/m$ is the square of the $(n_\lambda/2)$th elementary divisor of $S^\lambda$. But since $n!/n_\lambda$ is the 
product of all hook lengths of $\lambda$, and since $\lambda$ is symmetric, we conclude that $H/m$ is a square in $\Q$.
\end{Proposition}

\begin{Conjecture}
\label{ConjSym4}
The quotient $H/m$ is an integer.

\rm
Here is the list of the symmetric non-hook partitions whose elementary divisors are known so far {\bf\cite{L}}. By (\ref{CorTr3}), it is sufficient to list the 
respective first half of the elementary divisors. Thus, the last elementary divisor we give is the $(n_\lambda/2 + 1)$st, together with its multiplicity.

\begin{footnotesize}
\[
\begin{array}{|l|l|l|l|}\hline
\lambda     & m  & H  & \mb{elementary divisors} \\\hline
(2,2)       & 3  & 3  & 2^1\cdot 6^1 \\ 
(3,2,1)     & 5  & 5  & 1^4\cdot 3^4\cdot 15^4\cdots \\
(3^2,2)     & 15 & 15 & 8^{21}\cdot 120^{21} \\
(3^3)       & 15 & 15 & 24^{21}\cdot 360^{21} \\
(4,2,1^2)   & 7  & 7  & 2^{14}\cdot 8^{31}\cdot 56^{31}\cdots \\
(4,3,2,1)   & 21 & 21 & 1^{41}\cdot 3^{176}\cdot 15^{167}\cdot 315^{167}\cdots \\
(4,3^2,1)   & 21 & 21 & 8^{372}\cdot 40^{222}\cdot 840^{222}\cdots \\
(4^2,2^2)   & 35 & 35 & 8^{131}\cdot 24^{353}\cdot 72^{836}\cdot 2520^{836}\cdots\\
(4^2,3,2)   & 35 & 35 & 16^{428}\cdot 48^{2081}\cdot 144^{1781}\cdot 5040^{1781}\cdots\\
(5,2,1^3)   & 9  & 9  & 6^{56}\cdot 30^{168}\cdot 270^{168}\cdots \\
(5,3,2,1^2) & 3  & 27 & 2^{1046}\cdot 4^{460}\cdot 8^{44}\cdot 16^{386}\cdot 48^{1560}\cdot 144^{354}\cdot 432^{354}\cdots \\
(5,3^2,1^2) & 3  & 27 & 4^{2329}\cdot 20^{181}\cdot 40^{857}\cdot 120^{2559}\cdot 360^{2082}\cdot 1080^{2082}\cdots \\
(6,2,1^4)   & 11 & 11 & 24^{208}\cdot 48^2\cdot 144^{840}\cdot 1584^{840}\cdots \\
(7,2,1^5)   & 13 & 13 & 120^{792}\cdot 840^{3960}\cdot 10920^{3960}\cdots \\\hline
\end{array}
\]
\end{footnotesize}
\end{Conjecture}

Let $\alpha$ be the product of the strictly upper diagonal hook lengths of $\lambda$.
Let $\gamma$ be the first elementary divisor of the Gram matrix of $S^\lambda$. We denote $E := \End_{\sZ\Al_n} S^\lambda$ and $E_{\sQ} := \Q\ts_\sZ E$.
Symmetry of $\lambda$ gives
\[
S^\lambda\;\iso_{\,\sZ\Al_n}\; S^{\lambda,-} \;\auf{\mb{\scr\bf\cite[\rm 6.7]{J}}}{\iso}_{\hspace*{-2mm}\sZ\Sl_n}\; S^{\lambda',\ast}
\; =\; S^{\lambda,\ast}\;\iso\; S^{\lambda,\#} ,
\]
and so 
\[
(S^\lambda\lraa{\beta} S^\lambda)\; :=\; (S^\lambda\hra S^{\lambda,\#}\lraiso S^\lambda)\;\in\; E
\]
(cf.\ section \ref{SubSecUnimod}). By definition of $\beta$, we have $\rho^\lambda(\sigma)\beta = \beta\rho^\lambda(\sigma)\eps_\sigma\,$, where 
$\rho^\lambda(\sigma)$ is the operation of $\sigma\in\Sl_n$ on $S^\lambda$.

\begin{Proposition}
\label{PropSym1}
Suppose $(-1)^{(n-s)/2}H$ not to be a square in $\Z$ and denote $\chi :=$ \linebreak[4] $\sqrt{(-1)^{(n-s)/2} H}$. 
We have an isomorphism $E_{\sQ} \iso \Q(\chi)$ (since a splitting field for $E_\sQ$ contains the ordinary character values given by 
{\bf\cite[\rm 2.5.12, 2.15.13]{JK}}), which we fix and use as an identification. Let $h\in\Z_{>0}$ be maximal such that $h^2$ divides $H$, and let 
$\chi_0 := \sqrt{(-1)^{(n-s)/2} H/h^2}$. The following hold.
\begin{itemize}
\item[(i)] We have $\gamma\; |\; \alpha h$.
\item[(ii)] If the ring of algebraic integers in $E_\sQ$ is given by $\Z[\chi_0]$, then $E = \Z[\alpha\gamma^{-1}\chi]$. If the ring of algebraic integers in 
$E_\sQ$ is given by $\Z[(1 + \chi_0)/2]$, then the following holds. If $\alpha\gamma^{-1}h$ is even, then $E = \Z[\alpha\gamma^{-1}\chi]$; if 
$\alpha\gamma^{-1}h$ is odd, then $E = \Z[\alpha\gamma^{-1}h(1 + \chi_0)/2]$ or $E = \Z[\alpha\gamma^{-1}\chi]$.
\item[(iii)] The elementary divisors of the Gram matrix of $S^\lambda$ and the elementary divisors of the operation of $\alpha\chi\in E$ on the $E$-module 
$S^\lambda$ coincide.
\end{itemize}

\rm
Using the map $\Sl_n\lra\Gal(\E_{\sQ}/\Q)$ given by conjugation by means of $\rho^\lambda$, we conclude that $\beta$ has trace zero, i.e.\ that 
$\beta = \pm\w\alpha\chi$ for some $\w\alpha\in\Q_{> 0}$. 

On the other hand, $\beta^2 = (-1)^{(n-s)/2}\w\alpha^2 H$ shows that the $i$th and the $(n_\lambda + 1 - i)$th elementary divisor of $S^\lambda$ multiply to 
give $\w\alpha^2 H$ for $i\in [1,n_\lambda]$. But by (\ref{CorTr3}), they multiply to give $n!/n_\lambda$, which is the product of the hook lengths of 
$\lambda$ {\bf\cite[\rm 20.1]{J}}. Therefore, $\w\alpha = \alpha$, which proves (iii).

Suppose $\Z[x]$ to be the ring of algebraic integers in a quadratic number field, and suppose $F\tm \Z[x]$ to be a subring of $\Z$-rank $2$. Choosing a 
$\Z$-linear basis $( u + vx, w + yx )$ of $F$, where $u,v,w,y\in\Z$, the existence of $1\in F$ shows that there are $a,b\in\Z$ such that $au + bw = 1$, 
$av + by = 0$. Base change by $\smatzz{a}{b}{-w}{u}$ yields a $\Z$-linear basis of the form $( 1, dx)$, $d\geq 1$, i.e.\ $E = \Z[dx]$.

Now, $\gamma^{-1}\alpha\chi$ is still contained in $E$. More precisely, $\alpha$ is minimal with this property, for otherwise $\gamma^{-1}\alpha\chi$ would not 
have $1$ as first elementary divisor.

Therefore, if the ring of algebraic integers in $E_\sQ$ is given by $\Z[\chi_0]$, minimality of $\alpha$ yields $d = \alpha\gamma^{-1}h$. If the ring of 
algebraic integers in $E_\sQ$ is given by $\Z[(1 + \chi_0)/2]$, minimality of $\alpha$ yields $d = \alpha\gamma^{-1}h$ if $d$ is odd, and 
$d = 2\alpha\gamma^{-1}h$ if $d$ is even. This proves (ii). Moreover, in all three cases $\alpha\gamma^{-1}h$ is integral, as asserted in (i). 
\end{Proposition}

\begin{Proposition}
\label{PropSym2}
Suppose $(-1)^{\frac{n-s}{2}}H$ to be a square in $\Z$, so in particular $n\con_4 s$. We have an isomorphism $E_{\sQ}\iso \Q\ti\Q$ (since the central-primitive
idempotents of $\C\Al_n$ belonging to the summands of $S_\sC^\lambda|_{\Al_n}$ already lie in $\Q\Al_n$, as the character values given in 
{\bf\cite[\rm 2.5.12, 2.15.13]{JK}} show), which we fix and use as an identification. The following hold.

\begin{itemize}
\item[(i)] We have $\gamma\; |\; \alpha\sqrt{H} \; (= \sqrt{n!/n_\lambda})$.
\item[(ii)] We have $E = \{ (a,b)\in \Z\ti\Z\; |\; a\con_d b\}$ with $d = 2\alpha\gamma^{-1}\sqrt{H}$ if $\alpha\gamma^{-1}$ is even, and either 
$d = 2\alpha\gamma^{-1}\sqrt{H}$ or $d = \alpha\gamma^{-1}\sqrt{H}$ if $\alpha\gamma^{-1}$ is odd.
\item[(iii)] The elementary divisors of the Gram matrix of $S^\lambda$ and those of the operation of $\alpha \sqrt{H}(1,-1)$ on the $E$-module $S^\lambda$ 
coincide.
\end{itemize}

\rm
The nontrivial automorphism of $E_\sQ$ turns $\beta$ into $-\beta$, whence $\beta = \pm\w\alpha (1,-1)$ for some $\w\alpha\in\Q_{>0}$. Since 
$\beta^2 = \w\alpha^2$, comparison with (\ref{CorTr3}) yields $\w\alpha = \alpha\sqrt{H}$, which proves (iii).

We have $E = \{ (a,b)\in \Z\ti\Z\; |\; a\con_d b\}$ for some $d\geq 1$. But $\alpha\gamma^{-1}\sqrt{H}$ is minimal in $\Q_{> 0}$ with 
$\alpha\gamma^{-1}\sqrt{H}(1,-1)\in E$, that is, $\alpha\gamma^{-1}\sqrt{H}$ is integral -- as claimed in (i) -- and minimal with 
$2\alpha\gamma^{-1}\sqrt{H}\con_d 0$. Thus, $d = \alpha\gamma^{-1}\sqrt{H}$ if $d$ is odd, and $d = 2\alpha\gamma^{-1}\sqrt{H}$ if $d$ is even, whence (ii).
\end{Proposition}

\pagebreak[4]

\begin{footnotesize}
\begin{quote}
\begin{Example}
\label{ExH2_5}\rm
Here are some examples of the subring $E\tm E_\sQ$.

\begin{center}
\begin{tabular}{|l|l|l|l|}\hline
$\lambda$   & $E$                                    & $\alpha$              & $\gamma$  \\\hline
$(3,2,1)$   & $\Z[3(1+\sqrt{5})/2]$                  & $3$                   & $1$         \\
$(4,1^3)$   & $\Z[\sqrt{-7}]$                        & $6$                   & $6$         \\
$(4,2,1^2)$ & $\Z[4\sqrt{-7}]$                       & $8$                   & $2$        \\
$(3^2,2)$   & $\Z[\sqrt{-15}]$                       & $8$                   & $8$         \\
$(3^3)$     & $\Z[\sqrt{-15}]$                       & $24$                  & $24$        \\
$(5,1^4)$   & $\{(a, b)\in\Z\ti\Z\; |\; a\con_6 b\}$ & $24$ ($\sqrt{H} = 3$) & $24$ \\\hline
\end{tabular}
\end{center}

According to (\ref{PropSym1} iii, \ref{PropSym2} iii), the main problem that remains to be solved is to determine the structure of $S^\lambda$ as an $E$-module. 
For example, the elementary divisors of $S^{(3,2,1)}$ are given by $1^4\cdot 3^4\cdot 15^4\cdot 45^4$, whereas on the free $E$-module $E^8$, multiplication by 
$\alpha\chi = 3\sqrt{5}$ has the elementary divisors $1^8\cdot 45^8$. Thus by (\ref{PropSym1} iii), $S^\lambda$ is not free over $E$, and not locally free at 
$3$ either.
\end{Example}
\end{quote}
\end{footnotesize}

\begin{Example}
\label{ExM3}
\rm
The asymmetric example $S^\mu$ of smallest dimension for which two successive elementary divisors have quotient bigger than the outer hook length 
$\mu_1 + \mu'_1 - 1$, is given by $\mu = (9,2^2,1)$. In fact, we obtain the elementary divisors 
$4^{792}\cdot 8^{3421}\cdot 120^{493}\cdot 960^{714}\cdot 2880^{586}$ {\bf\cite{L}}, whereas the outer hook length is $12$.
\end{Example}

The following remark on partitions of rectangular shape pertains in particular to partitions of quadratic shape.

\begin{Remark}
\label{RemM4}
Let $\mu$ be a partition such that $\mu_1 = \mu_h$, where $h = \mu'_1$. Let $\nu$ be the partition obtained from $\mu$ by removing the lower right node in 
its diagram, i.e.\ $\nu_h := \mu_h - 1$, $\nu_j := \mu_j$ for $j\neq h$. Then the elementary divisors of $S^\mu$ are given by the elementary divisors of $S^\nu$,
multiplied by the constant factor $h$.

\rm
We have $S^\mu|_{\Sl_{n-1}} \llaiso S^\nu$ {\bf\cite[\rm 9.3]{J}}, by sending a $\nu$-polytabloid to the $\mu$-polytabloid obtained by adding the entry $n$ in 
the lower right corner. Thus the invariant bilinear form on $S^\nu$ induces an invariant bilinear form on $S^\mu$, whence the Gram matrix of $S^\mu$ is a scalar 
multiple of the Gram matrix of $S^\nu$. This scalar is calculated to be $h$ by the Branching Theorem for Determinants {\bf\cite[\rm p.\ 225]{JM}}.
\end{Remark}

\subsection{A conjectural comparison of kernels}

\begin{quote}
\begin{footnotesize}
There seems to be a connection between elementary divisors and modular morphisms. Apart from the fact that both yield necessary conditions on the shape of the 
quasiblock (cf.\ section \ref{SubSecQuasi}), we do not know of any a priori reason for such a connection to exist.
\end{footnotesize}
\end{quote}

Let $h\geq 1$, let $k\in [1,h]$. Abbreviate $\nu(k) := (2^{h-k},1^{n - 2h + k})$ and $b := n - 2h + 1$. Let $F^{\nu(0)}$ be the free $\Z$-module on the 
$\nu(0)$-tableaux, endowed with the natural operation of $\Sl_n$, yielding $F^{\nu(0)}\iso\Z\Sl_n$. Consider the morphism
\[
\begin{array}{rcl}
F^{\nu(0)} & \lraa{f_k} & S^{\nu(k)}/(b+k)S^{\nu(k)} \\
{[a]}      & \lramaps   & \sum_{\zeta\in Z(a,k)} (-1)^{\Sigma(\zeta)} \spi{a^\zeta}\; , \\
\end{array}
\]
where $Z(a,k)$ is the set of subtuples of length $k$ of the second column of $[a]$, where $[a^\zeta]$ is the $\mu$-tableau given by removing $\zeta$ from the 
second column of $[a]$ and appending it at the bottom of the first column, and where $\Sigma(\zeta)$ is the sum of the positions of the entries of $\zeta$. So 
for example, if the second column is given by $(3,5,6,7,9)$, and if $\zeta = (5,7,9)$, then $\Sigma(\zeta) = 2 + 4 + 5$.

Let $p$ be a prime. Given $m\geq 1$, we denote $\ilog_p(m) := \max\{ i\in\Z_{\geq 0}\; |\; p^i \leq m\}$. Let $m? := \prod_\mb{\scr $p$ prime} p^{\silog_p(m)}$ 
(`$m$ quasi-factorial', e.g.\ $4? = 12$). From {\bf\cite[\rm 4.9]{K00}} we take the morphism 
\[
\begin{array}{rcl}
S^{\nu(0)} & \lraa{s_k} & S^{\nu(k)}/(b+k)S^{\nu(k)} \\
\spi{a}    & \lramaps   & k?\cdot ([a]f_k)\; . \\
\end{array}
\]
We will form an intersection of the kernels of a slight variation of these morphisms -- the factor $k?$ can sometimes be lowered on a smaller domain.

Let $S^{\nu(0),\,\cap\, 0} := S^{\nu(0)}$, and $\w s_1 := s_1$.

Assume submodules $S^{\nu(0),\,\cap\, l-1}\tm S^{\nu(0)}$ and morphisms $S^{\nu(0),\,\cap\, l-1}\lraa{\w s_l} S^{\nu(l)}/(b+l) S^{\nu(l)}$ to be constructed for 
$l\in [1,k-1]$.

Let $S^{\nu(0),\,\cap\, k-1} := \Cap_{l\in [1,k-1]} \Kern\w s_l$. Let
\[
\begin{array}{rcl}
S^{\nu(0),\,\cap\, k-1}           & \lraa{\w s_k} & S^{\nu(k)}/(b+k) S^{\nu(k)} \\
\sum_{\spi{a}}t_{\spi{a}}\spi{a}  & \lramaps      & \frac{k?}{\gcd(k?,\gamma_{k-1})} \sum_{\spi{a}}t_{\spi{a}} ([a] f_k)\; ,
\end{array}
\]
where the $t_{\spi{a}}$ are coefficients in $\Z$, and where $\gamma_{k-1}$ is the first elementary divisor of $S^{\nu(0),\,\cap\, k-1}\hra S^{\nu(0)}$.

\begin{Conjecture}
\label{ConM5}
The kernel of
\[
S^{\nu(0)}\;\;\lraa{\eta}\;\; S^{\nu(0),\ast}/{\text\frac{h!\, (n-h+1)!}{n-2h+1}}\,S^{\nu(0),\ast} 
= S^{\nu(0),\ast}/{\text\frac{n!}{\rks_\ssZ S^{\nu(0)}}}\,S^{\nu(0),\ast}
\]
is given by the intersection
\[
S^{\nu(0),\,\cap\, h}\; =\; \Cap_{k\in [1,h]} \Kern \w s_k\; .
\]
\end{Conjecture}

Conjecture (\ref{ConM5}) holds if $n\leq 8$ (direct computation), or if $h = 1$ (\ref{ThHook}). If $h = 2$, the kernel is contained in the intersection 
(\ref{Th2^2}).

\subsection{A final remark on quasiblocks}
\label{SubSecQuasi}

Let $p$ be a prime, let $R = \Z_{(p)}$, and let $\eps^\lambda$ denote the central-primitive idempotent of $\Q\Sl_n$ belonging to $S_\sQ^\lambda$, for $\lambda$ 
a partition of $n$. The quasiblocks $\eps^\lambda R\Sl_n$ of $R\Sl_n$ may be viewed as the `building blocks' of the locally integral representation theory
of the symmetric group. The Gram matrix determines a certain part of the structure of such a quasiblock, sometimes it even describes it in its entirety. 

Wedderburn's isomorphism $\Q\Sl_n\lraiso\prod_{\lambda}\End_\sQ S_\sQ^\lambda$ restricts to $R\Sl_n\hra\prod_{\lambda} \End_R S_R^\lambda\,$, which in turn 
induces $\eps^\lambda R\Sl_n\hraa{\rho^\lambda}\End_R S_R^\lambda\,$, given by sending $\eps^\lambda\sigma$ to the operation $\rho^\lambda(\sigma)$ of 
$\sigma\in\Sl_n$ on $S_R^\lambda$. A necessary condition for an $R$-linear endomorphism of $S_R^\lambda$ to lie in the image of this embedding is its 
compatibility with the $R\Sl_n$-linear embedding $S_R^\lambda\hraa{\eta} S_R^{\lambda,\ast}$. Namely, if $S_R^\lambda\lraa{\phi}S_R^\lambda$ is in the image 
of $\rho^\lambda$, then there exists an $R$-linear operation $S_R^{\lambda,\ast}\lraa{\psi}S_R^{\lambda,\ast}$ such that $\phi\eta = \eta\psi$.
In terms of matrices, this means that the image $\rho^\lambda(R\Sl_n) = \rho^\lambda(\eps^\lambda R\Sl_n)$ is contained in 
$\Gamma_R^\lambda := (R)_{n_\lambda}\cap \left(G^\lambda (R)_{n_\lambda} (G^{\lambda})^{-1}\right)$, where $G^\lambda$ denotes the Gram matrix.
For instance, we obtain 
\[
\eps^{(2,1)} \Z_{(3)}\Sl_3 \;\;\lraisoa{\rho^{(2,1)}}\;\; \Gamma_{\sZ_{(3)}}^{(2,1)}\; =\;
(\Z_{(3)})_2\cap \left(\smatzz{1}{0}{0}{3}(\Z_{(3)})_2\smatzz{1}{0}{0}{3}^{-1}\right) \; =\; \smatzz{\sZ_{(3)}}{\sZ_{(3)}}{(3)}{\sZ_{(3)}}\; .
\]
But in general, the inclusion 
\[
\eps^\lambda R\Sl_n\;\hraa{\rho^\lambda}\; \Gamma_R^\lambda
\]
is not surjective -- as an example we may take $\lambda = (3,2,1)$ over $R = \Z_{(3)}$.

If all indecomposable projective $\eps^\lambda R\Sl_n$-lattices are simple, that is, if the quasiblock is a tiled order, we do not know whether
\begin{itemize}
\item[(i)] $\;\Jac(\eps^\lambda R\Sl_n)\;\hraa{\rho^\lambda}\; \Jac(\Gamma_R^\lambda)\;$,
\end{itemize}
and if so, whether
\begin{itemize}
\item[(ii)] there exists a $K\geq 0$ such that for all $k\geq K$
\[
\Jac^k(\eps^\lambda R\Sl_n)\;\lraisoa{\rho^\lambda}\; \Jac^k(\Gamma_R^\lambda)\; .
\]
\end{itemize}
Cf.\ {\bf\cite[\rm 6.1.26]{K99}}. Several examples of tiled quasiblocks are given in {\bf\cite{P80,P83}}. An example of a tiled quasiblock that requires $K = 2$ 
is given in {\bf\cite[\rm 5.6.12 ii]{N}}. 

\parskip0.0ex
\begin{footnotesize}

\parskip1.2ex

\vss

M.\ K\"unzer, G.\ Nebe\\
Abt.\ Reine Mathematik\\
Universit\"at Ulm\\
D-89069 Ulm\\
kuenzer@mathematik.uni-ulm.de \\
nebe@mathematik.uni-ulm.de \\
\end{footnotesize}

\renewcommand{\thefootnote}{\fnsymbol{footnote}}
\footnotetext[0]{{\bf Note added in proof.} Matthew Fayers independently obtained (\ref{CorT4}) and (\ref{PropTr1}) using different methods. See 
{\sc M.\ Fayers,} {\it On the structure of Specht modules,} J.\ Lond.\ Math.\ Soc.\ (2) 67, p.\ 85--102, 2003.} %
\footnotetext[0]{{\bf Note added in proof.} {\sc G.\ E.\ Murphy} obtained (\ref{CorT4}) in an explicit form that avoids matrix inversion 
(Thesis, Queen Mary College, London, 1987). It is stated in loc.\ cit.\ (6.27) as being dependent of the conjecture loc.\ cit.\ (6.24), which in turn has been 
proven by Schaper {\bf\cite[\rm p.\ 60]{Sch}}.} %
\renewcommand{\thefootnote}{\arabic{footnote}}


\end{document}